\documentclass[notheoremnums]{dpreprint}

\newtheorem{theorem}{Theorem}

\newtheorem{conjecture}[theorem]{Conjecture}
\newtheorem{lemma}{Lemma}[section]

\newtheorem{proposition}[lemma]{Proposition}

\theoremstyle{definition}
\newtheorem{definition}[lemma]{Definition}

\newtheorem{remark}[lemma]{Remark}
\newtheorem{notation}[lemma]{Notation}

\newcommand{\len}[1]{|#1|}

\begin{document}

\dtitle[Word complexity cutoffs for mixing properties]{Word complexity cutoffs for mixing properties of subshifts}
\dauthor{Darren Creutz}{darren.a.creutz@vanderbilt.edu}{Vanderbilt University}{}
\datewritten{\today}

\keywords{Symbolic dynamics, word complexity, mixing, partial mixing}
\subjclass{Primary: 37B10; Secondary 37A25}

\dabstract{
Investigating word complexity cutoffs for partial and light mixing in subshifts, we construct partially mixing rank-one subshifts whose word complexity functions $p(q)$ satisfy
$
\liminf_{q\to\infty}\nicefrac{p(q)}q<2+\epsilon
$
for any prescribed $\epsilon>0$ and, conversely, we prove that every lightly mixing rank-one subshift satisfies
$
\liminf_{q\to\infty}\nicefrac{p(q)}q>2,
$
establishing $2$ as the optimal limiting threshold in the rank-one setting. More quantitatively, we show that every $\alpha$-partially mixing rank-one subshift satisfies
$
\liminf_{q\to\infty}\nicefrac{p(q)}q\ge 2+\nicefrac{\alpha}{(1+\alpha)}.
$
Our constructions are a combination of staircase and Chacon-like rank-one stages; we introduce the independently interesting notion of mixing spans and meshes to establish partial mixing while retaining low word complexity.
}

\makepreprint

\dsectionstar{Introduction}

Subshifts are a natural and well-studied class of transformations: transformations on closed subsets $X \subseteq \mathcal{A}^{\mathbb{Z}}$ which are invariant under the shift map (here $\mathcal{A}$, the alphabet, is an arbitrary finite set). A key aspect of subshifts is that there is a finer notion of complexity than for general transformations, namely we can define the \textbf{word complexity} function $p(q) = |\{ w~\text{is a subword of some $x \in X$}~:~|w| = q \}|$ where $|w|$ is the length. Entropy, the standard measure of complexity for general transformations, cannot distinguish qualitative behavior, such as various forms of measure-theoretic mixing, among zero entropy transformations; however, word complexity can. Questions of what the minimal possible word complexity for subshifts having mixing properties date back to Ferenczi in the 1990s, e.g.~\cite{ferenczi1996rank},\cite{ferenczichacon} and are an active area of current research.

A \textbf{word complexity cutoff} is a quantitative statement about the asymptotic nature of $p(q)$ for which some specified qualitative property of the subshift can occur. For example, in \cite{creutzmixing}, the author showed that the word complexity cutoff for measure-theoretically strong mixing is precisely that $p(q)$ is superlinear: any strongly mixing subshift has $\lim p(q)/q = \infty$ and for any $f(q) \to \infty$ there exists a strongly mixing subshift for which $p(q) < q f(q)$ (for sufficiently large $q$). For another example, the related question for weak mixing was resolved by the author and R.~Pavlov \cite{creutzpavlov1.5} who showed that $\limsup p(q)/q < 1.5$ implies discrete spectrum, hence precludes weak mixing, and \cite{creutzpavlovds}
where the existence of a weakly mixing subshift with $\limsup \nicefrac{p(q)}{q} = 1.5$ and $\liminf \nicefrac{p(q)}{q} = 1$ was established.

A striking feature of recent work on word complexity is the appearance of a sequence of thresholds of the form $(k+1)/k$, in several different senses, at which increasingly restrictive structural and dynamical properties cease to be forced. In particular, $3/2$ is the sharp threshold for the occurrence of weak mixing in terms of $\limsup p(q)/q$ \cite{creutzpavlov1.5}, while $4/3$ marks a transition in the substitutive structures available to low-complexity subshifts \cite{creutzpavlovds}. Further restrictions on substitution incidence matrices arise below $5/4$, see \cite{creutzpavlovds}, specifically the discussion following Proposition 3.1. At the other end of this hierarchy, strong mixing requires superlinear word complexity, i.e.~``$1/0$". The results of this paper identify $2 = 2/1$ as the corresponding limiting threshold for partial mixing within the rank-one class and provide evidence for the analogous threshold for light mixing in general.

A transformation $T$ is \textbf{$\alpha$-partially mixing} for $0 < \alpha < 1$ if for all measurable sets $A$ and $B$,
\[
\liminf \mu(T^{n}A \cap B) \geq \alpha \mu(A) \mu(B)
\]
A transformation $T$ is \textbf{lightly mixing} if for all positive measure sets $A$ and $B$,
$
\liminf \mu(T^{n}A \cap B) > 0
$. Strong mixing is defined as $\lim \mu(T^{n}A \cap B) = \mu(A)\mu(B)$ and weak mixing as the existence of a sequence of times along which that limit holds; as is well-known the implications
\[
\text{strong mixing} \quad \Rightarrow \quad \text{partial mixing} \quad \Rightarrow \quad \text{light mixing} \quad \Rightarrow \quad \text{weak mixing}
\]
hold and none of the reverse implications hold.

\subsection*{Main results}

Here we establish a quantitative, in terms of word complexity, result regarding partial mixing as it fits into the general word complexity cutoff framework:

\begin{theorem}\label{examples}
For every $\epsilon > 0$, there exists a partially mixing subshift such that $\liminf \nicefrac{p(q)}{q} < 2 + \epsilon$.
\end{theorem}

While weak mixing is possible with $\liminf p(q)/q = 1$ and strong mixing forces $\liminf p(q)/q = \infty$, our work here shows that light mixing can occur with $\liminf p(q)/q < 2+\epsilon$. Building on this, we establish that, in the context of rank-one subshifts, $2$ is the optimal bound (and in fact cannot be attained by a rank-one subshift though we suspect it can be attained by a non-rank-one lightly mixing subshift):


\begin{theorem}\label{R1nolight}
For any lightly mixing rank-one subshift, $\liminf \nicefrac{p(q)}{q} > 2$.
\end{theorem}

We also establish lower bounds for word complexity of arbitrary $\alpha$-partially mixing rank-one subshifts demonstrating that stronger partial mixing forces increasing lower bounds on word complexity.

\begin{theorem}\label{alphaoptimalR1}
For any rank-one subshift which is $\alpha$-partially mixing, $\liminf \nicefrac{p(q)}{q} \geq 2 + \nicefrac{\alpha}{(\alpha + 1)}$.
\end{theorem}

\subsection*{Further directions}

Motivated by these results, we propose the following:

\begin{conjecture}\label{conj0}
For any subshift which is $\alpha$-partially mixing, $\liminf \nicefrac{p(q)}{q} \geq 2 + \nicefrac{\alpha}{(\alpha+1)}$.
\end{conjecture}

\begin{conjecture}\label{conj1}
For any lightly mixing subshift, $\liminf \nicefrac{p(q)}{q} \geq 2$.
\end{conjecture}

The apparent discrepancy with the above and our result is intentional, and in fact we suspect that in the rank-two setting, equality can be attained just as equality with $1.5$ can be attained in the weakly mixing rank-two setting despite the fact that weakly mixing rank-one subshifts cannot attain the $1.5$ bound.

\begin{conjecture}
There exists a lightly mixing subshift with $\liminf \nicefrac{p(q)}{q} = 2$.
\end{conjecture}

Our lightly mixing examples with $\liminf p(q)/q < 2 + \epsilon$ necessarily have $\limsup p(q)/q = \infty$ due to the nature of their construction; we conjecture that such is unavoidable in general when constructing lightly mixing subshifts:

\begin{conjecture}\label{conj2}
For any lightly mixing subshift, $\limsup \nicefrac{p(q)}{q} = \infty$.
\end{conjecture}

%

\dsubsectionstar{Mixing spans and meshes}

A key new ingredient in the proof of Theorem \ref{examples} is the notion of mixing spans and meshes, which is of independent interest and worth highlighting. A mesh is, in essence, a way of taking specific subsets of Rokhlin towers which themselves are organized into levels and therefore retain much of the usefulness of towers when estimating mixing-type asymptotics while being significantly more general.

Our examples essentially depend on meshes in that they are built as a combination of staircase-type constructions and Chacon-type constructions; the staircase stages introduce some amount of mixing and meshes allow us to show that said mixing survives in a suitable sense through the Chacon-like stages of the construction. In terms of word complexity, the staircase stages, as used in our constructions, unavoidably increase the complexity in a superlinear fashion and the Chacon-like stages are required to then reduce the word complexity arbitrarily close to the prescribed bounds; as such the notion of meshes is crucial to the argument.

We introduce the notion of (uniform) mixing stages of rank-one transformations which involves spans of times during which a fixed set is mixed uniformly into every mesh laid over every sufficiently large column.  These new ingredients, combined with  staircase mixing techniques going back to \cite{adamsfriedman} and used, for example, in \cite{adams1998smorodinsky}, \cite{CreutzSilva2004}, \cite{CreutzSilva2010}, \cite{CreutzStochasticStaircases2021}, \cite{creutzmixing}, form the bulk of the proof of partial mixing (the other ingredient being the nonrigidity of our constructions).

\subsection*{Organization of the paper}

Section 1 establishes precise definitions and recalls the standard theorems on subshifts and rank-one constructions, Section 2 defines the class of examples we will use to establish Theorem \ref{examples} and Section 3 provides the word complexity estimates for these examples. Sections 4 and 5 establish partial mixing, with Section 4 being largely devoted to developing the theory of meshes as they relate to mixing estimates and Section 5 being more focused on our specific class of examples. Section 6 concretely exhibits the examples proving Theorem \ref{examples} and the paper concludes with Section 7 where the lower bounds in Theorems \ref{R1nolight} and \ref{alphaoptimalR1} are demonstrated.

\subsection*{Acknowledgments}

The author thanks Terry Adams profusely for an extremely enjoyable multi-day conversation during which, among other things, he suggested the author consider this problem and made multiple useful observations about how one might prove partial mixing for subshifts while retaining something resembling low word complexity.

\section{Definitions and preliminaries}

\begin{definition}
A \textbf{subshift} on the finite set $\mathcal{A}$ is any subset $X \subset \mathcal{A}^{\mathbb{Z}}$ which is closed in the product topology and shift-invariant: for all $x = (x_{n})_{n \in \mathbb{Z}} \in X$ and $k \in \mathbb{Z}$, the translate $(x_{n+k})_{n \in \mathbb{Z}}$ of $x$ by $k$ is also in $X$.
\end{definition} 

\begin{definition}
A \textbf{word} is any element of $\mathcal{A}^\ell$ for some $\ell$, the \textbf{length} of $w$, written $\len{w}$.  A word $w$ is a \textbf{subword} of a word or bi-infinite sequence $x$ if there exists $k$ so that $w_{i} = x_{i+k}$ for all $1 \leq i \leq \len{w}$.  A word $u$ is a \textbf{prefix} of $w$ when $u_{i} = w_{i}$ for $1 \leq i \leq \len{u}$ and a word $v$ is a suffix of $w$ when $v_{i} = w_{i+\len{w}-\len{v}}$ for $1 \leq i \leq \len{v}$.
\end{definition}

For words $v,w$, we denote by $vw$ their concatenation--the word obtained by following $v$ immediately by $w$. We write such concatenations with product or exponential notation, e.g.~$\prod_i w_i$ or $0^n$.

\begin{definition}
The \textbf{language} of a subshift $X$ is $\mathcal{L}(X) = \{ w : \text{$w$ is a subword of some $x \in X$} \}$.
\end{definition}

\begin{definition}
The \textbf{word complexity function} of a subshift $X$ over $\mathcal{A}$ is the function $p_X: \mathbb{N} \rightarrow \mathbb{N}$ defined by $p_X(q) = |\mathcal{L}(X) \cap \mathcal{A}^q|$, the number of words of length $q$ in the language of $X$.
\end{definition}

When $X$ is clear from context, we suppress the subscript and just write $p(n)$.

For subshifts on the alphabet $\{0,1\}$, we consider:

\begin{definition}
The set of \textbf{right-special} words is
$
\mathcal{L}^{RS}(X) = \{ w \in \mathcal{L}(X) : w0, w1 \in \mathcal{L}(X) \}
$.
\end{definition}

Cassaigne \cite{cassaigne} showed the well-known: $
p(q) = p(m) + \sum_{\ell = m}^{q-1} |\{ w \in \mathcal{L}^{RS} : \len{w} = \ell \}|
$
for $m < q$.

\begin{definition}
A \textbf{transformation} $T$ is a measurable map on a standard Borel or Lebesgue measure space $(Y,\mathcal{B},\mu)$ that is measure-preserving: $\mu(T^{-1}B) = \mu(B)$ for all $B \in \mathcal{B}$.
\end{definition}

\begin{definition}Two transformations $T$ on $(Y,\mathcal{B},\mu)$ and $T^{\prime}$ on $(Y^{\prime},\mathcal{B}^{\prime},\mu^{\prime})$ are \textbf{measure-theoretically isomorphic} when there exists 
a bijective map $\phi$ between full measure subsets $Y_0 \subset Y$ and $Y'_0 \subset Y'$ where 
$\mu(\phi^{-1} A) = \mu'(A)$ for all measurable $A \subset Y'_0$ and $(\phi \circ T)(y) = (T' \circ \phi)(y)$ for all $y \in Y_0$.
\end{definition}

%
%

\begin{definition}
A transformation $T$ on a finite measure space is \textbf{lightly mixing} if for all positive measure sets $A$ and $B$,
\[
\liminf \mu(T^{n}A \cap B) > 0
\]
\end{definition}

\begin{definition}
A transformation $T$ on a finite measure space is \textbf{$\alpha$-partially mixing} for $0 < \alpha < 1$ if for all measurable sets $A$ and $B$,
\[
\liminf \mu(T^{n}A \cap B) \geq \alpha \mu(A) \mu(B)
\]
\end{definition}

A \textbf{rank-one transformation} is a transformation $T$ constructed by ``cutting and stacking".  Here $Y$ represents a (possibly infinite) interval, $\mathcal{B}$ is the induced $\sigma$-algebra from $\mathbb{R}$, and $\mu$ is Lebesgue measure. We give a brief description, referring the reader to 
\cite{fghsw21} or  \cite{silva2008invitation} for more details.

The transformation is defined inductively on larger and larger portions of the space through Rohlin towers or \textbf{columns}, denoted $C_n$. Each column $C_n$ consists of \textbf{levels} $I_{n,j}$ where $0 \leq j < h_{n}$ is the height of the level within the column. All levels $I_{n,j}$ in $C_n$ are 
intervals with the same length, $\mu(I_{ n})$, and the total number of levels in a column is the \textbf{height} of the column, denoted by $h_n$. The transformation $T$ is defined on all levels $I_{n,j}$ except the top one $I_{n, h_n - 1}$ by sending each $I_{n,j}$ to $I_{n,j+1}$ using the unique order-preserving affine map.

Start with $C_1=[0,1)$ with height $h_1=1$. To obtain $C_{n+1}$ from $C_n$, we require a \textbf{cut sequence}, $\{r_n\}$ such that $r_n \geq 1$ for all $n$. Make $r_n$ vertical cuts of $C_n$ to create $r_n+1$ \textbf{subcolumns} of equal width. Denote a \textbf{sublevel} of $C_n$ by $I_{n,j}^{[i]}$ where $0 \leq a < h_{n}$ is the height of the level within that column, and 
$i$ represents the position of the subcolumn, where $i=0$ represents the leftmost subcolumn and $i=r_n$ is the rightmost subcolumn. After cutting $C_n$ into subcolumns, add extra intervals called \textbf{spacers} on top of each subcolumn to function as levels of the next column. The \textbf{spacer sequence}, $\{s_{n,i}\}$ such that $0 \leq i \leq r_{n}$ and $s_{n,i} \geq 0$, specifies how many sublevels to add above each subcolumn. Spacers are the same width as the sublevels, act as new levels in the column $C_{n+1}$, and are taken to be the leftmost intervals in $[1,\infty)$ not in $C_{n}$.
After the spacers are added, stack the subcolumns with their spacers right on top of left, i.e.~so that $I_{n,0}^{[i+1]}$ is directly above $I_{n,h_{n}-1}^{[i]}$. This gives the next column, $C_{n+1}$.

Each column $C_n$ defines $T$ on $\bigcup_{j = 0}^{h_n - 2} I_{n,j}$ and 
the partially defined map $T$ on $C_{n+1}$ agrees with that of $C_n$, extending the definition of $T$ to a portion of the top level of $C_n$ where it was previously undefined. Continuing this process gives the sequence of columns $\{C_1, \dots, C_n, C_{n+1}, \dots\}$ and $T$ is then the limit of the partially defined maps. 

Though this construction could result in $Y$ being an infinite interval with infinite Lebesgue measure, $Y$ has finite measure if and only if 
$\sum_{n} \frac{1}{(r_{n}+1)h_{n}}\sum_{i=0}^{r_{n}} s_{n,i} < \infty$, see \cite{CreutzSilva2010}.
All rank-one transformations we define satisfy this condition, and for convenience we renormalize so that $Y = [0,1)$. Every rank-one transformation is ergodic and invertible.

The reader should be aware that we are making $r_{n}$ cuts and obtaining $r_{n}+1$ subcolumns (following Ferenczi \cite{ferenczi1996rank}), while other papers (e.g.~\cite{CreutzStochasticStaircases2021}) use $r_{n}$ as the number of subcolumns.

For a rank-one transformation defined as above, we define a subshift $X(T)$ on the alphabet $\{ 0,1 \}$ which is measure-theoretically isomorphic to $T$ provided $T$ is on a finite measure space:

\begin{definition}
The \textbf{symbolic model} $X(T)$ of a rank-one transformation $T$ is given by the sequence of words: $B_0 = 0$ and 
\begin{equation*}\label{rk1word}
B_{n+1}=B_n1^{s_{n,0}}B_n1^{s_{n,1}}\cdots B_{n} 1^{s_{n,r_n}}= \prod_{i=0}^{r_n} B_n1^{s_{n,i}}
\end{equation*}
and $X(T)$ is the set of all bi-infinite sequences such that every subword is a subword of some $B_n$.
\end{definition}

The words $B_n$ are a symbolic coding of the column $C_n$: $0$ represents $C_1$ and $1$ represents the spacers.  There is a natural measure associated to $X(T)$:

\begin{definition}
The \textbf{empirical measure} for a symbolic model $X(T)$ of a rank-one transformation $T$ is the measure $\nu$ defined by, for each word $w$,
\[
\nu([w]) = \lim_{n \rightarrow \infty} \frac{|\{1 \leq j \leq \len{B_{n}} - \len{w} + 1 \ : \ {B_n}_{[j,j+\len{w})} = w\}|}{\len{B_{n}} - \len{w} + 1}
\]
\end{definition}

Danilenko \cite{danilenkopr} (combined with \cite{deljunco} and \cite{kalikow}) proved that the symbolic model $X(T)$ of a finite measure rank-one subshift, equipped with its empirical measure, is measure-theoretically isomorphic to the cut-and-stack construction (see \cite{adamsferenczipeterson17}; see \cite{fghsw21} for the full generality including odometers).

Due to this isomorphism, we move back and forth between rank-one and symbolic model terminology as needed and write $\mathcal{L}(T)$ for the language of $X(T)$.

\section{Staircase-Chacon-like hybrid transformations}

\begin{definition}
Given a word $B_{n}$ representing a column in a rank-one transformation, say that $B_{n+1}$ is of \textbf{staircase form} when
\[
B_{n+1} = \left(\prod_{j=0}^{r_{n}-1}B_{n}1^{c_{n}+j}\right)B_{n}
\]
for some $r_{n} \geq 2$ and $c_{n} \geq 0$.
\end{definition}

\begin{definition}
Given a word $B_{n}$ representing a column in a rank-one transformation, say that $B_{n+1}$ is
of \textbf{Chacon-like form} when 
\[
B_{n+1} = \left(\prod_{j=0}^{r_{n}-1}B_{n}1^{c_{n} + \iota_{n,j}M_{n}}\right)B_{n}
\]
where $\iota_{n,j} \in \mathbb{N} \cup \{ 0 \}$ for all $n$ and $r_{n} \geq 2$ and $c_{n} \geq 0$ and $M_{n} \geq 1$.
\end{definition}

\begin{definition}
Let $\{ \mathrm{Form}_{n} \}$ be a sequence of forms all either staircase or Chacon-like and let $\{ M_{n} \}$ be a sequence of positive integers for those $n$ for which $B_{n+1}$ is of Chacon-like form and $\{ r_{n} \}$ a sequence of positive integers.

Set $B_{0} = 0$ and $c_{0} = 1$ and $r_{0} = 3$.  Given, $B_{n}$, construct $B_{n+1}$ to be of the specified form.

If $B_{n+1}$ and $B_{n+2}$ are of staircase form, set $c_{n+1} = c_{n} + r_{n}$.

If $B_{n+1}$ is of staircase form and $B_{n+2}$ is of Chacon-like form, set $c_{n+1} \geq c_{n} + r_{n}$ minimal such that $M_{n+1}$ divides $|B_{n+1}| + c_{n+1}$.

If $B_{n+1}$ is of Chacon-like form and $B_{n+2}$ is of staircase form, let $c_{n+1} \geq |B_{\ell}| + 1 + M_{n+1}$ be any integer
where $\ell < n$ is the largest stage that was of staircase form.

If $B_{n+1}$ and $B_{n+2}$ are of Chacon-like form, set $c_{n+1} = c_{n}$.

The resulting transformation is an \textbf{elevated staircase-Chacon-like hybrid transformation}.
\end{definition}

\begin{notation}
For an elevated staircase-Chacon-like hybrid transformation,
\[
\tilde{h}_{n} = |B_{n}| + c_{n}
\]
\end{notation}

\begin{lemma}\label{heights}
Let $T$ be an elevated staircase-Chacon-like hybrid transformation.  If $B_{n+1}$ is constructed using the staircase operation then
\[
\tilde{h}_{n+1} = (r_{n}+1)\tilde{h}_{n} + \frac{1}{2}r_{n}(r_{n}-1) + c_{n+1} - c_{n}
\]
and if $B_{n+1}$ is constructed using a Chacon-like operation then
\[
\tilde{h}_{n+1} = (r_{n}+1)\tilde{h}_{n} + s_{n}M_{n} + c_{n+1} - c_{n}
\]
where $s_{n}$ is the number of subcolumns over which non-elevation spacers are placed.
\end{lemma}
\begin{proof}
In general, if $s_{n,j}$ is the number of non-elevation spacers added above each subcolumn then $\tilde{h}_{n+1} = (r_{n}+1)h_{n} + r_{n}c_{n} + \sum_{j=0}^{r_{n}-1} s_{n,j} + c_{n+1}$ since no spacers (elevating or otherwise) are added above the rightmost subcolumn.

For $n$ such that $B_{n+1}$ is constructed using the staircase operation, $s_{n,j} = j$ so $\tilde{h}_{n+1} = (r_{n}+1)\tilde{h}_{n} + \frac{1}{2}r_{n}(r_{n}-1) + c_{n+1} - c_{n}$.

For $n$ such that $B_{n+1}$ is constructed using a Chacon-like operation, $\sum s_{n,j} = s_{n}M_{n}$ so the claim follows similarly.
\end{proof}

\begin{definition}
Say that $B_{n+1}$ is constructed using the $\gamma$-Chacon-like form for $\gamma \geq 2$ if
\[
B_{n+1} = \left(B_{n}1^{c_{n}}1^{M_{n}}\right)^{\gamma}B_{n}1^{c_{n}}\left(B_{n}1^{c_{n}}1^{M_{n}}\right)^{\gamma}B_{n}
\]
for some $M_{n} \geq 1$.
\end{definition}

\section{Word complexity estimates}

\begin{proposition}\label{3prime}
Let $T$ be an elevated staircase-Chacon-like hybrid transformation and $n$ be such that $B_{n+1}$ is of staircase form.  
Then $p(|B_{n+1}|+c_{n+1})$ is completely determined by the cutting and stacking operations up to $B_{n+1}$.
\end{proposition}
\begin{proof}
Let $w$ be a right-special word with $|w| \leq |B_{n+1}| + c_{n+1}$.  Since $1^{|w|}$ will appear regardless, we may assume $w$ has $01^{z}$ as a suffix for some $z \geq 1$ (as $00$ is not in the language because $c_{n} \geq 1$).

If $z < c_{n+1}$ then $w$ must be a subword of $1^{c_{n}}B_{n+1}$ because $1^{c_{n+1}}v$ is not right-special for any proper prefix $v$ of $B_{n+1}$ (as $1^{c_{n+1}}$ is always followed by $B_{n+1}$).

If $z > c_{n+1}$ then, since $1^{c_{n+1}+M_{n+1}}$ is not right-special and $c_{N} \geq |B_{n+1}|$ for the first $N > n+1$ for which $B_{N+1}$ is constructed using the staircase operation, then $w = 1^{|w|}$.

If $z = c_{n+1}$ then $w$ has $01^{c_{n+1}}$ as a suffix hence shares a suffix with $1^{c_{n+1}}B_{n+1}1^{c_{n+1}}$ and is therefore a suffix of that word.

Therefore all right-special words of length at most $|B_{n+1}| + c_{n+1}$ are determined by $B_{n+1}$ and $c_{n+1}$.
\end{proof}

\begin{proposition}\label{00}
Let $w$ be a right-special word in the language of an elevated staircase-Chacon-like hybrid transformation such that $|B_{n}| < |w|$ for $n$ such that $B_{n+1}$ is of Chacon-like form and such that $|w| \leq C_{N}$ for the smallest $N > n$ such that $B_{N}$ is of staircase form.
Then one of the following holds:
\begin{itemize}
\item $w = 1^{|w|}$
\item $w$ shares a suffix with $1^{c_{n}}B_{n}1^{c_{n}}$
\end{itemize}
\end{proposition}
\begin{proof}
If $w = 1^{|w|}$ then we are in the first case so assume that $w$ has $01^{z}$ as a suffix for some $z \geq 0$.  By construction, $00 \notin \mathcal{L}$ so $z \geq 1$.

Let $m$ such that $c_{m} \leq z < c_{m+1}$.
Here $m$ is unique by the construction of $\{ c_{n} \}$.  This follows as when $B_{m+1}$ is of staircase form, $c_{m}$ is strictly larger than any spacer amounts added at any prior stage and $c_{m+1}$ is strictly greater than any spacer amounts added in $B_{m+1}$ (and any values outside the ranges between $c_{n}$ and $c_{n+1}$ can never be the suffix of a right-special word).

Suppose $m \geq N$.  Then $|w| \leq c_{N} \leq c_{m}$ so $w = 1^{|w|}$.  Therefore $m < N$.

Suppose $m < n$ and $B_{m+1}$ is constructed using the staircase operation.
Then $w$ ends properly inside a copy of $B_{m+1}$.  Then $w$ has a suffix of the form $u1^{c_{m+1}}v$ for some word $u$ and some prefix $v$ of $B_{m}$ and $v$ has $1^{c_{m}+j}$ as a suffix for some $0 \leq j < r_{m}$.  This follows as $|w| > |B_{n}|$.  Then $w$ has $1^{c_{m+1}}B_{m}1^{c_{m}}B_{m}1^{c_{m}+1}\ldots B_{m}1^{c_{m}+j-1} B_{m}1^{c_{m}+j}$ as a suffix.  But that word is never followed by a $1$ so is not right-special.  Therefore $m < n$ is impossible in this case.

Likewise, if $m < n$ and $m$ is prior to a staircase operation $t$ which is prior to $n$ then $w$ must contain $1^{c_{t+1}}B_{t}1^{c_{t}}B_{t}1^{c_{t}+1}\ldots B_{t}1^{c_{t}+j-1} B_{t}1^{c_{t}+j}$ followed by a prefix of $B_{t}$.  But that word is never right-special.

Therefore $m$ and $n$ are not separated by a staircase operation.  
So we may assume that $c_{n} \leq z < c_{N}$.  Therefore $z$ is either $c_{n}$ or $c_{n} + M_{n}$ by construction.  Since $01^{c_{n} + M_{n} + 1}$ is not in the language (as $C_{N} > c_{n} + M_{n} + 1$), then $z = c_{n}$.  So $w$ has $01^{c_{n}}$ as a suffix.  Then by construction $w$ shares a suffix with $B_{n}1^{c_{n}}$ hence with $1^{c_{n}}B_{n}1^{c_{n}}$.
\end{proof}

\begin{lemma}\label{c1}
Let $w$ be a right-special word in the language of an elevated staircase-Chacon-like hybrid transformation such that $|B_{n}| < |w|$ for $n$ such that all of $B_{n}$ and $B_{n+1}$ and $B_{n+2}$ are of $\gamma$-Chacon-like form and such that $w$ shares a suffix with $11^{c_{n}}B_{n}1^{c_{n}}$. 
Then one of the following holds:
\begin{itemize}
\item $w$ shares a suffix with $B_{n+1}1^{c_{n}}$
\item $w$ is a suffix of $\left(B_{n}1^{c_{n}}1^{M_{n}}\right)^{2\gamma})B_{n}1^{c_{n}}$
\end{itemize}
\end{lemma}
\begin{proof}
Here $w$ shares a suffix with $1^{c_{n}}B_{n}1^{c_{n}+M_{n}}B_{n}1^{c_{n}}$ since $11^{c_{n}}B_{n}1^{c_{n}}$ only appears as a suffix of that word.
By construction, $c_{n} = c_{n+1}$ and $M_{n} = M_{n+1}$.

Suppose $w$ shares a suffix with $01^{c_{n}}B_{n}1^{c_{n}+M_{n}}B_{n}1^{c_{n}}$.
Then $w$ and $1^{c_{n}}B_{n}1^{c_{n}}B_{n}1^{c_{n}+M_{n}}B_{n}1^{c_{n}}$ have a common suffix but $1^{c_{n}}B_{n}1^{c_{n}}B_{n}1^{c_{n}+M_{n}}B_{n}1^{c_{n}}B_{n}$ is not in the language (since $\gamma > 1$) so this would mean $w$ is not right-special.

So $w$ shares a suffix with $11^{c_{n}}B_{n}1^{c_{n}+M_{n}}B_{n}1^{c_{n}}$ hence with $1^{c_{n}+M_{n}}B_{n}1^{c_{n}+M_{n}}B_{n}1^{c_{n}}$.

If $w$ shares a suffix with $01^{c_{n}}(B_{n}1^{c_{n}+M_{n}})^{t}B_{n}1^{c_{n}+M_{n}}B_{n}1^{c_{n}}$ for any $t \geq 1$ then $w0$ shares a suffix with $01^{c_{n}}(B_{n}1^{c_{n}+M_{n}})^{t}B_{n}1^{c_{n}+M_{n}}B_{n}1^{c_{n}}0$ which only appears as a subword of $(B_{n}1^{c_{n}})^{2}(B_{n}1^{c_{n}+M_{n}})^{t+1}(B_{n}1^{c_{n}})^{2}$ meaning that $t \geq \gamma$.  Therefore $w$ shares a suffix with $(B_{n}1^{c_{n}+M_{n}})^{\gamma}B_{n}1^{c_{n}}$.

Consider when $w$ shares a suffix with $11^{c_{n}}(B_{n}1^{c_{n}+M_{n}})^{\gamma}B_{n}1^{c_{n}+M_{n}}B_{n}1^{c_{n}}$.  That word only appears as a subword of $B_{n+1}1^{c_{n}+M_{n}}B_{n+1}1^{c_{n}}$ and is always preceded by $(B_{n}1^{c_{n}+M_{n}})^{\gamma}B_{n}1^{M_{n}-1}$.  Therefore $w$ shares a suffix with $(B_{n}1^{c_{n}+M_{n}})^{2\gamma+1}B_{1}1^{c_{n}}$.

Since $(B_{n}1^{c_{n}+M_{n}})^{2\gamma+1}B_{1}1^{c_{n}}$ is not right-special (because $(B_{n}1^{c_{n}+M_{n}})^{2\gamma+2}$ is not in the language), then $w$ is a suffix of $(B_{n}1^{c_{n}+M_{n}})^{2\gamma}B_{1}1^{c_{n}}$ and we are in the second claimed case.

Now consider the case when $w$ shares a suffix with $01^{c_{n}}(B_{n}1^{c_{n}+M_{n}})^{\gamma}B_{n}1^{c_{n}}$.  Then $w$ shares a suffix with $B_{n}1^{c_{n}}(B_{n}1^{c_{n}+M_{n}})^{\gamma}B_{n}1^{c_{n}}$ which is always preceded by $(B_{n}1^{c_{n}+M_{n}})^{\gamma}B_{n}1^{c_{n}}$, in which case we are in the second case claimed, or by preceded $1^{c_{n}+M_{n}}$ in which case $w$ shares a suffix with $((B_{n}1^{c_{n}+M_{n}})^{\gamma}B_{n}1^{c_{n}})^{2}$ so with $B_{n+1}1^{c_{n}}$ and we are in the first claimed case.
\end{proof}

\begin{lemma}\label{c2}
Let $w$ be a right-special word in the language of an elevated staircase-Chacon-like hybrid transformation such that $|B_{n}| < |w|$ for $n$ such that all of $B_{n}$ and $B_{n+1}$ and $B_{n+2}$ are of $\gamma$-Chacon-like form and such that $w$ shares a suffix with $01^{c_{n}}B_{n}1^{c_{n}}$.

Then
$w$ is a suffix of $(B_{n-1}1^{c_{n-1}}1^{M_{n-1}})^{\gamma-1}B_{n-1}1^{c_{n-1}}B_{n}1^{c_{n}}$.
\end{lemma}
\begin{proof}
Here $w$ shares a suffix with $1^{c_{n}}B_{n}1^{c_{n}}B_{n}1^{c_{n}}$. 
By construction, $c_{n-1} = c_{n}$ and $M_{n-1} = M_{n}$.

Since $B_{n}1^{c_{n}}B_{n}$ only appears as a subword of $1^{c_{n}+M_{n}}B_{n}1^{c_{n}}B_{n}1^{c_{n}+M_{n}}$, then
 $1^{c_{n}}B_{n}1^{c_{n}}B_{n}1^{c_{n}}B_{n}$ is not in the language so $1^{c_{n}}B_{n}1^{c_{n}}B_{n}1^{c_{n}}0$ is not so $w$ must be a proper suffix of $1^{c_{n}}B_{n}1^{c_{n}}B_{n}1^{c_{n}}$.

Consider when $w$ shares a suffix with $01^{c_{n}}(B_{n-1}1^{c_{n}+M_{n}})^{\gamma}B_{n-1}1^{c_{n}}B_{n}1^{c_{n}}$.  That word only appears when the leading zero is the tail zero in the first $(B_{n-1}1^{c_{n}+M_{n}})^{\gamma}B_{n-1}$ which is in the first $B_{n}$ of $B_{n}1^{c_{n}}B_{n}1^{c_{n}}$.  Then $01^{c_{n}}(B_{n-1}1^{c_{n}+M_{n}})^{\gamma}B_{n-1}1^{c_{n}}B_{n}1^{c_{n}}0$ is not in the language so $w$ is a suffix of $1^{c_{n}}(B_{n-1}1^{c_{n}+M_{n}})^{\gamma-1}B_{n-1}1^{c_{n}}B_{n}1^{c_{n}}$ as claimed.

Now consider when $w$ shares a suffix with $11^{c_{n}}(B_{n-1}1^{c_{n}+M_{n}})^{\gamma}B_{n-1}1^{c_{n}}B_{n}1^{c_{n}}$.  Then it shares a suffix with $1^{c_{n}+M_{n}}(B_{n-1}1^{c_{n}+M_{n}})^{\gamma}B_{n-1}1^{c_{n}}B_{n}1^{c_{n}}$ so, as $1^{c_{n}+M_{n}}$ is always preceded by $B_{n-1}$ (as it is a suffix of $B_{n}$), then $w$ shares a suffix with $1^{c_{n}}B_{n-1}1^{c_{n}+M_{n}}(B_{n-1}1^{c_{n}+M_{n}})^{\gamma}B_{n-1}1^{c_{n}}B_{n}1^{c_{n}}$ but $(B_{n-1}1^{c_{n}+M_{n}})^{\gamma+1}$ is not in the language so this is impossible.
\end{proof}

\begin{proposition}\label{1full}
Let $w$ be a right-special word in the language of an elevated staircase-Chacon-like hybrid transformation such that $|B_{n}| < |w|$ for $n$ such that all of $B_{n}$ and $B_{n+1}$ and $B_{n+2}$ are of $\gamma$-Chacon-like form and such that $|w| \leq C_{N}$ for the smallest $N > n$ such that $B_{N}$ is of staircase form.
Then one of the following holds:
\begin{itemize}
\item $w = 1^{|w|}$
\item $w$ shares a suffix with $B_{n+1}1^{c_{n}}$
\item $w$ is a suffix of $\left(B_{n}1^{c_{n}}1^{M_{n}}\right)^{2\gamma}B_{n}1^{c_{n}}$
\item $w$ is a suffix of $(B_{n-1}1^{c_{n-1}}1^{M_{n-1}})^{\gamma-1}B_{n-1}1^{c_{n-1}}B_{n}1^{c_{n}}$
\end{itemize}
\end{proposition}
\begin{proof}
By Proposition \ref{00}, if $w \ne 1^{|w|}$ then $w$ shares a suffix with $1^{c_{n}}B_{n}1^{c_{n}}$.  By Lemmas \ref{c1} and \ref{c2} then the claim follows.
\end{proof}

\begin{proposition}\label{1full2}
Let $T$ be an elevated staircase-Chacon-like hybrid transformation such that all of $B_{n}$ and $B_{n+1}$ and $B_{n+2}$ are of $\gamma$-Chacon-like form.  Then
\[
p(q+1) - p(q) \leq \left\{ \begin{array}{ll}
3 \quad &\text{if $\tilde{h}_{n} < q \leq \frac{3}{2}\tilde{h}_{n}$} \\[0.1em]
2 \quad &\text{if $\frac{3}{2}\tilde{h}_{n} < q \leq \frac{1}{2}\tilde{h}_{n+1}$} \\[0.1em]
3 \quad &\text{if $\frac{1}{2}\tilde{h}_{n+1} < q \leq \tilde{h}_{n+1} - \tilde{h}_{n}$} \\[0.1em]
2 \quad &\text{if $\tilde{h}_{n+1} - \tilde{h}_{n} < q \leq \tilde{h}_{n+1}$}
\end{array} \right.
\]
\end{proposition}
\begin{proof}
Let $\tilde{h}_{n} < q \leq \tilde{h}_{n+1}$.  The words $1^{q}$ and the suffix of $B_{n+1}1^{c_{n}}$ of length $q$ are both right-special.  By Proposition \ref{1full}, the only remaining right-special words are suffixes of $(B_{n}1^{c_{n}+M_{n}})^{2\gamma}B_{n}1^{c_{n}}$ and $(B_{n-1}1^{c_{n}+M_{n}})^{\gamma-1}B_{n-1}1^{c_{n}}B_{n}1^{c_{n}}$.

Since $(B_{n}1^{c_{n}+M_{n}})^{2\gamma}B_{n}1^{c_{n}}$ and $B_{n+1}1^{c_{n}}$ share a common suffix of length $(\gamma + 1)\tilde{h}_{n} + \gamma M_{n} = \frac{1}{2}\tilde{h}_{n+1}$, the only new right-special words from that word have length $\frac{1}{2}\tilde{h}_{n+1} < q \leq (2\gamma + 1)\tilde{h}_{n} + 2\gamma M_{n} = \tilde{h}_{n+1} - \tilde{h}_{n}$.

Since $(B_{n-1}1^{c_{n}+M_{n}})^{\gamma-1}B_{n-1}1^{c_{n}}B_{n}1^{c_{n}}$ shares a suffix of length $\tilde{h}_{n} + c_{n}$ with $B_{n+1}1^{c_{n}}$, the only new right-special words from that word have length $\tilde{h}_{n} < q \leq \gamma \tilde{h}_{n-1} + \gamma M_{n} + \tilde{h}_{n} = \frac{1}{2}\tilde{h}_{n}$.
\end{proof}


\begin{proposition}\label{2full}
Let $T$ be an elevated staircase-Chacon-like hybrid such that the $\gamma$-Chacon-like operation is used to construct $B_{m}$ through $B_{n+2}$ for $n > m$.  Then
\[
\frac{p(\frac{1}{2}\tilde{h}_{n+1})}{\frac{1}{2}\tilde{h}_{n+1}} <
2 + \frac{1}{\gamma + 1} + 3^{m-n} \frac{p(\tilde{h}_{m})}{\tilde{h}_{m}}
\]
\end{proposition}
\begin{proof}
Since $B_{n}$ and $B_{n+1}$ and $B_{n+2}$ are all constructed using the $\gamma$-Chacon-like operation, by Proposition \ref{1full2},
\[
p(\tilde{h}_{n+1}) - p(\tilde{h}_{n})
\leq 2(\tilde{h}_{n+1} - \tilde{h}_{n}) + \frac{1}{2}\tilde{h}_{n} + \frac{1}{2}\tilde{h}_{n} + \frac{1}{2}\tilde{h}_{n+1} - \tilde{h}_{n}
= \frac{5}{2}\tilde{h}_{n+1} - \frac{5}{2}\tilde{h}_{n}
\]
and therefore, by the telescoping nature of the series,
\[
p(\tilde{h}_{n}) \leq \frac{5}{2}\tilde{h}_{n} - \frac{5}{2}\tilde{h}_{m} + p(\tilde{h}_{m})
\]
Again by Proposition \ref{1full2},
\[
p\left(\frac{1}{2}\tilde{h}_{n+1}\right) - p(\tilde{h}_{n})
\leq 2\left(\frac{1}{2}\tilde{h}_{n+1} - \tilde{h}_{n}\right) + \frac{1}{2}\tilde{h}_{n}
\]
and therefore
\[
p\left(\frac{1}{2}\tilde{h}_{n+1}\right) \leq
\tilde{h}_{n+1} - 2\tilde{h}_{n} + \frac{1}{2}\tilde{h}_{n} + \frac{5}{2}\tilde{h}_{n} - \frac{5}{2}\tilde{h}_{m} + p(\tilde{h}_{m})
\leq \tilde{h}_{n+1} + \tilde{h}_{n} + p(\tilde{h}_{m})
\]
meaning that
\[
\frac{p(\frac{1}{2}\tilde{h}_{n+1})}{\frac{1}{2}\tilde{h}_{n+1}}
\leq 2 + \frac{\tilde{h}_{n}}{(\gamma + 1)\tilde{h}_{n} + \gamma M_{n}} + \frac{p(\tilde{h}_{m})}{\frac{1}{2}\tilde{h}_{n+1}}
< 2 + \frac{1}{\gamma + 1} + \frac{\tilde{h}_{m}}{\tilde{h}_{n}} \frac{p(\tilde{h}_{m})}{\tilde{h}_{m}}
\]
and the claim follows since $\tilde{h}_{t+1} > 3\tilde{h}_{t}$ (as $\gamma \geq 2$) for all $t$.
\end{proof}

\section{Mixing spans and meshes}

For a column, $C_{n}$, recall that $I_{n,0}$ is the base of the column and $I_{n,j} = T^{j}I_{n,0}$ for $0 \leq j < \tilde{h}_{n}$ are the levels (including the elevating levels).

For a level $I_{n,j}$ and $1 \leq i \leq r_{n}$, recall that $I_{n,j}^{[i]}$ is the $i^{th}$ sublevel, i.e.~the $\nicefrac{1}{r_{n}+1}$ fraction of $I_{n,j}$ in the $i^{th}$ subcolumn.

\begin{notation}
For a measurable set $B$, let $\chi_{B} = \bbone_{B} - \mu(B)$.

For $J_{n,0} \subseteq I_{n,0}$, let $J_{n,0}^{[i]} := J_{n,0} \cap I_{n,0}^{[i]}$ for $0 \leq i \leq r_{n}$.
\end{notation}

\begin{definition}
Let $I_{n,0}$ be the base level in a tower $C_{n}$ of a rank-one transformation.   If $J_{n,0} \subseteq I_{n,0}$ such that for all $0 \leq i \leq r_{n}$, it holds that $J_{n,0} \cap I_{n,0}^{[i]}$ has the property that $J_{n,0}^{[i]} = T^{i\tilde{h}_{n} + \sum_{z=0}^{i-1}s_{n,z}}J_{n,0}^{[0]}$ then $J_{n,0}$ is a \textbf{mesh prototype} for $I_{n,0}$.

Given a mesh prototype $J_{n,0}$ for $I_{n,0}$, for $0 \leq j < \tilde{h}_{n}$ and $0 \leq i \leq r_{n}$, write
\[
J_{n,j} = T^{j}J_{n,0} \quad\quad\text{and}\quad\quad J_{n,j}^{[i]} = T^{j}\left(J_{n,0} \cap I_{n,0}^{[i}\right)
\]
and say that $\bigsqcup_{j=0}^{\tilde{h}_{n}-1} J_{n,j}$ is a \textbf{mesh} for $C_{n}$ and that $J_{n,j}$ are the \textbf{mesh levels}.
\end{definition}

Note that $C_{n}$ is a mesh for itself by putting $J_{n,j} = I_{n,j}$ and also that the intersection of two meshes for $C_{n}$ is again a mesh for $C_{n}$.

\begin{definition}
For a mesh $\bigsqcup_{j=0}^{\tilde{h}_{n}-1}J_{n,j}$ of some column $C_{n}$ of a rank-one transformation, the \textbf{size} of the mesh is
\[
\frac{\mu\left(\bigsqcup_{j=0}^{\tilde{h}_{n}-1}J_{n,j}\right)}{\mu(C_{n})}
= \frac{\mu(J_{n,0})}{\mu(I_{n,0})}
\]
\end{definition}

\begin{lemma}\label{mesh1}
Let $J_{n,0}$ be a mesh prototype for some base $I_{n,0}$ of some column $C_{n}$ of a rank-one transformation.  Let $A$ be a union of levels in $C_{n}$.  Then for $0 \leq j < \tilde{h}_{n}$ and $0 \leq i \leq r_{n}$,
\[
\mu(J_{n,j}^{[i]} \cap A) = \frac{1}{r_{n}+1}\mu(J_{n,j} \cap A)
\]
\end{lemma}
\begin{proof}
Since $A$ is a union of levels in $C_{n}$, either $I_{n,j} \subseteq A$ or $I_{n,j} \cap A = \varnothing$.
Since $J_{n,j}^{[i]} \subseteq I_{n,j}^{[i]} \subseteq I_{n,j}$, if $I_{n,j} \cap A = \varnothing$ then $J_{n,j} \cap A = \varnothing$ so $\mu(J_{n,j}^{[i]} \cap A) = 0 = \frac{1}{r_{n}+1}0 = \frac{1}{r_{n}+1}\mu(J_{n,j} \cap A)$.  Likewise, if $I_{n,j} \subseteq A$ then $J_{n,j} \subseteq A$ and since $\mu(J_{n,j}^{[i]}) = \mu(J_{n,0}^{[0]})$ (as $T$ is measure-preserving and from the definition of a mesh) is constant over all $0 \leq i \leq r_{n}$ and $j$ (for a fixed $n$), then $\mu(J_{n,j}^{[i]} \cap A) = \mu(J_{n,j}^{[i]}) = \frac{1}{r_{n}+1}\mu(J_{n,j}) = \frac{1}{r_{n}+1}\mu(J_{n,j} \cap A)$.
\end{proof}

\begin{lemma}\label{meshdrop}
Let $\bigsqcup J_{n+1,j}$ be a mesh of $C_{n+1}$ for some column of a rank-one transformation.  Then $C_{n} \cap \bigsqcup J_{n+1,j}$ is a mesh of $C_{n}$.
\end{lemma}
\begin{proof}
Let $s_{n,i}$ for $1 \leq i \leq r_{n}$ be the number of non-elevation spacers added above each subcolumn of $C_{n}$.  Then for $0 \leq j < \tilde{h}_{n}$,
\[
I_{n,j} = \bigsqcup_{i=0}^{r_{n}} I_{n,j}^{[i]} = \bigsqcup_{i=0}^{r_{n}} I_{n+1,j+i\tilde{h}_{n} + \sum_{t=0}^{i-1} s_{n,t}}
\]
and set
\[
J_{n,j} := 
I_{n,j} \cap \bigsqcup_{m=0}^{\tilde{h}_{n+1}-1} J_{n+1,m}
= \bigsqcup_{i=0}^{r_{n}} \left(I_{n+1,j+i\tilde{h}_{n} + \sum_{t=0}^{i-1} s_{n,t}} \cap J_{n+1,j+i\tilde{h}_{n} + \sum_{t=0}^{i-1} s_{n,t}}\right)
\]
Then
\[
J_{n,j}^{[i]} = I_{n+1,j+i\tilde{h}_{n} + \sum_{t=0}^{i-1} s_{n,t}} \cap J_{n+1,j+i\tilde{h}_{n} + \sum_{t=0}^{i-1} s_{n,t}}
= T^{i\tilde{h}_{n} + \sum_{0=1}^{i-1} s_{n,t}} J_{n,j}^{[0]}
\]
meaning $\bigsqcup J_{n,j}$ form a mesh of $C_{n}$.
\end{proof}

\begin{lemma}\label{ts}
Let $T$ be a rank-one transformation and $A$ be a measurable set and let $\epsilon > 0$ and $K \in \mathbb{N}$.  Let $\bigsqcup_{j=0}^{\tilde{h}_{n}-1}J_{n,j}$ be a mesh for $C_{n}$ for some $n$.  If for
all $1 \leq k \leq K$,
\[
\sum_{j=0}^{\tilde{h}_{n}-1}\left|\mu(T^{k\tilde{h}_{n}}J_{n,j} \cap A) - \mu J_{n,j} \mu A \right| < \epsilon
\]
then for all $\tilde{h}_{n} \leq t \leq K\tilde{h}_{n}$,
\[
\sum_{j=0}^{\tilde{h}_{n}-1}\left|\mu(T^{t}J_{n,j} \cap A) - \mu J_{n,j} \mu A \right| < 2\epsilon
\]
\end{lemma}
\begin{proof}
Let $\tilde{h}_{n} \leq t < K\tilde{h}_{n}$.  Write $t = k\tilde{h}_{n} + z$ for the unique $0 \leq z < \tilde{h}_{n}$ and $1 \leq k < K$.  Then
\begin{align*}
\sum_{j=0}^{\tilde{h}_{n}-1}&\left|\mu(T^{t}J_{n,j} \cap A) - \mu J_{n,j} \mu A \right| \\
&=
\sum_{j=\tilde{h}_{n}-z}^{\tilde{h}_{n}-1}\left|\mu(T^{t}J_{n,j} \cap A) - \mu J_{n,j} \mu A \right| +
\sum_{j=0}^{\tilde{h}_{n}-z-1}\left|\mu(T^{t}J_{n,j} \cap A) - \mu J_{n,j} \mu A \right| \\
&= 
\sum_{j=0}^{z-1}\left|\mu(T^{(k+1)\tilde{h}_{n}}J_{n,j} \cap A) - \mu J_{n,j} \mu A \right| +
\sum_{j=z}^{\tilde{h}_{n}-1}\left|\mu(T^{k\tilde{h}_{n}}J_{n,j} \cap A) - \mu J_{n,j} \mu A \right|
< \epsilon + \epsilon \qedhere
\end{align*}
\end{proof}

\begin{lemma}\label{mix00}
Let $T$ be a rank-one  transformation and $n$ be such that $B_{n+1}$ is constructed using the staircase operation.  Let $A$ be a union of levels in $C_{n}$ and $\bigsqcup_{j=0}^{\tilde{h}_{n}-1}J_{n,j}$ be a mesh for $C_{n}$.  Then for $\tilde{h}_{n} > j \geq r_{n}$,
\[
\left|\mu(T^{\tilde{h}_{n}}J_{n,j} \cap A) - \mu J_{n,j} \mu A\right| \leq 
\int_{J_{n,j}}\left|\frac{1}{r_{n}+1}\sum_{i=0}^{r_{n}-1}\chi_{A} \circ T^{-i}\right| d\mu + \frac{1}{r_{n}+1}\mu J_{n,j}
\]
\end{lemma}
\begin{proof}
For $j \geq r_{n}$ and $0 \leq i < r_{n}$,
\[
T^{\tilde{h}_{n}}J_{n,j}^{[i]} = T^{-i}J_{n,j}^{[i+1]} = J_{n,j-i}^{[i+1]}
\]
and since $A$ is a union of levels, by Lemma \ref{mesh1},
\[
\mu(T^{\tilde{h}_{n}}J_{n,j}^{[i]} \cap A) = \mu(J_{n,j-i}^{[i+1]} \cap A) 
= \frac{1}{r_{n}+1}\mu(J_{n,j-i} \cap A) = \frac{1}{r_{n}+1}\mu(T^{-i}J_{n,j} \cap A)
\]
and therefore
\[
\mu(T^{\tilde{h}_{n}}J_{n,j} \cap A) = \frac{1}{r_{n}+1}\sum_{i=0}^{r_{n}-1}\mu(T^{-i}J_{n,j} \cap A) +
\mu(T^{\tilde{h}_{n}}J_{n,j}^{[r_{n}]} \cap A)
\]
Then, since $\mu(T^{\tilde{h}_{n}}J_{n,j}^{[r_{n}]} \cap A) \leq \frac{1}{r_{n}+1}\mu J_{n,j}$,
\begin{align*}
\left|\mu(T^{\tilde{h}_{n}}J_{n,j} \cap A) - \mu J_{n,j} \mu A\right| &\leq 
\left| \frac{1}{r_{n}+1}\sum_{i=0}^{r_{n}-1}\mu(T^{-i}J_{n,j} \cap A) - \mu J_{n,j} \mu A \right| + \frac{1}{r_{n}+1}\mu J_{n,j} \\
&= \left| \int_{J_{n,j}} \frac{1}{r_{n}+1} \sum_{i=0}^{r_{n}-1} \chi_{A} \circ T^{-i}~d\mu \right| + \frac{1}{r_{n}+1}\mu J_{n,j} \\
&\leq \int_{J_{n,j}}\left|\frac{1}{r_{n}+1}\sum_{i=0}^{r_{n}-1}\chi_{A} \circ T^{-i}\right| d\mu + \frac{1}{r_{n}+1}\mu J_{n,j} \qedhere
\end{align*}
\end{proof}

\begin{lemma}\label{mix00inv}
Let $T$ be a rank-one  transformation and $n$ be such that $B_{n+1}$ is constructed using the staircase operation.  Let $A$ be a union of levels in $C_{n}$ and $\tilde{C}_{n} = \bigsqcup_{j=0}^{\tilde{h}_{n}-1}J_{n,j}$ be a mesh for $C_{n}$.  Then for $\tilde{h}_{n} > j \geq r_{n}$,
\[
\left|\mu(T^{\tilde{h}_{n}}I_{n,j} \cap A \cap \tilde{C}_{n}) - \mu J_{n,j} \mu A \right| \leq 
\int_{J_{n,j}}\left|\frac{1}{r_{n}+1}\sum_{i=0}^{r_{n}-1}\chi_{A} \circ T^{-i}\right| d\mu + \mu (T^{\tilde{h}_{n}}I_{n,j}^{[r_{n}]} \cap A \cap \tilde{C}_{n})
\]
\end{lemma}
\begin{proof}
For $j \geq r_{n}$ and $0 \leq i < r_{n}$,
\[
T^{\tilde{h}_{n}}I_{n,j}^{[i]} = T^{-i}I_{n,j}^{[i+1]} = I_{n,j-i}^{[i+1]}
\]
since $I_{n,j-i}^{[i+1]} \cap A \cap \tilde{C}_{n} \subseteq J_{n,j-i}^{[i+1]}$ and since $A$ is a union of levels, by Lemma \ref{mesh1},
\[
\mu(T^{\tilde{h}_{n}}I_{n,j}^{[i]} \cap A \cap \tilde{C}_{n})
= \mu(I_{n,j-i}^{[i+1]} \cap A \cap \tilde{C}_{n})
= \mu(J_{n,j-i}^{[i+1]} \cap A)
= \frac{1}{r_{n}+1}\mu(J_{n,j-i} \cap A) = \frac{1}{r_{n}+1}\mu(T^{-i}J_{n,j} \cap A)
\]
and the claim follows as in Lemma \ref{mix00}.
\end{proof}

\begin{lemma}\label{mix0}
Let $T$ be a rank-one  transformation and $n$ be such that $B_{n+1}$ is constructed using the staircase operation.   Let $A$ be a union of levels in $C_{n}$ and $\tilde{C}_{n} = \bigsqcup_{j=0}^{\tilde{h}_{n}-1}J_{n,j}$ be a mesh for $C_{n}$.  Then
\[
\sum_{j=0}^{\tilde{h}_{n}-1}\left|\mu(T^{\tilde{h}_{n}}J_{n,j} \cap A) - \mu J_{n,j} \mu A \right|
\leq \int_{\tilde{C}_{n}} \left|\frac{1}{r_{n}+1}\sum_{i=0}^{r_{n}-1} \chi_{A} \circ T^{-i}\right| d\mu
+ \left(\frac{1}{r_{n}+1}
+ \frac{r_{n}}{\tilde{h}_{n}}\right)\mu \tilde{C}_{n} 
\]
and
\[
\sum_{j=0}^{\tilde{h}_{n}-1}\left|\mu(T^{\tilde{h}_{n}}I_{n,j} \cap A \cap \tilde{C}_{n}) - \mu J_{n,j} \mu A \right|
\leq \int_{\tilde{C}_{n}} \left|\frac{1}{r_{n}+1}\sum_{i=0}^{r_{n}-1} \chi_{A} \circ T^{-i}\right| d\mu
+ \sum_{j=0}^{\tilde{h}_{n}-1}\mu (T^{\tilde{h}_{n}}I_{n,j}^{[r_{n}]} \cap \tilde{C}_{n})
+ \frac{r_{n}}{\tilde{h}_{n}}\mu \tilde{C}_{n}
\]
\end{lemma}
\begin{proof}
Since the bottom $r_{n}$ levels in the mesh have measure at most $\frac{r_{n}}{\tilde{h}_{n}}\mu \tilde{C}_{n}$, by Lemma \ref{mix00},
\begin{align*}
\sum_{j=0}^{\tilde{h}_{n}-1}\left|\mu(T^{\tilde{h}_{n}}J_{n,j} \cap A) - \mu J_{n,j} \mu A \right|
&\leq \sum_{j=r_{n}}^{\tilde{h}_{n}-1} \left|\mu(T^{\tilde{h}_{n}}J_{n,j} \cap A) - \mu J_{n,j} \mu A\right| + \frac{r_{n}}{h_{n}}\mu\tilde{C}_{n} \\
&\leq \sum_{j=r_{n}}^{\tilde{h}_{n}-1} \int_{J_{n,j}}\left|\frac{1}{r_{n}+1}\sum_{i=0}^{r_{n}-1}\chi_{A}\circ  T^{-i}\right| d\mu + \frac{1}{r_{n}+1}\mu\tilde{C}_{n} + \frac{r_{n}^{2}}{\tilde{h}_{n}}\mu\tilde{C}_{n}  \\
&\leq \int_{\tilde{C}_{n}}\left|\frac{1}{r_{n}+1}\sum_{i=0}^{r_{n}-1}\chi_{A}\circ T^{-i}\right| d\mu + \frac{1}{r_{n}+1}\mu \tilde{C}_{n} + \frac{r_{n}}{\tilde{h}_{n}}\mu \tilde{C}_{n}
\end{align*}
where we used that the levels are disjoint to consolidate the integral.  Similarly, using Lemma \ref{mix00inv},
\begin{align*}
\sum_{j=0}^{\tilde{h}_{n}-1}&\left|\mu(T^{\tilde{h}_{n}}I_{n,j} \cap A \cap \tilde{C}_{n}) - \mu J_{n,j} \mu A \right|
\leq \sum_{j=r_{n}}^{\tilde{h}_{n}-1} \left|\mu(T^{\tilde{h}_{n}}I_{n,j} \cap A \cap \tilde{C}_{n}) - \mu J_{n,j} \mu A\right| + \frac{r_{n}}{\tilde{h}_{n}}\mu \tilde{C}_{n} \\
&\quad\quad\quad\quad\leq \sum_{j=r_{n}}^{\tilde{h}_{n}-1} \int_{J_{n,j}}\left|\frac{1}{r_{n}+1}\sum_{i=0}^{r_{n}-1}\chi_{A}\circ  T^{-i}\right| d\mu +  \sum_{j=0}^{\tilde{h}_{n}-1}\mu (T^{\tilde{h}_{n}}I_{n,j}^{[r_{n}]} \cap \tilde{C}_{n}) + \frac{r_{n}^{2}}{\tilde{h}_{n}}\mu \tilde{C}_{n}  \\
&\quad\quad\quad\quad\leq \int_{\tilde{C}_{n}}\left|\frac{1}{r_{n}+1}\sum_{i=0}^{r_{n}-1}\chi_{A}\circ T^{-i}\right| d\mu +  \sum_{j=0}^{\tilde{h}_{n}-1}\mu (T^{\tilde{h}_{n}}I_{n,j}^{[r_{n}]} \cap \tilde{C}_{n}) + \frac{r_{n}}{\tilde{h}_{n}}\mu \tilde{C}_{n} \qedhere
\end{align*}
\end{proof}

\begin{lemma}\label{mix0all}
Let $T$ be a rank-one  transformation and $n$ be such that $B_{n+1}$ is constructed using the staircase operation.  Let $A$ be a union of levels in $C_{n}$ and $\tilde{C}_{n}^{(q)}$ be a collection of disjoint meshes of $C_{n}$.  Then
\[
\sum_{q} \sum_{j=0}^{\tilde{h}_{n}-1}\left|\mu(T^{\tilde{h}_{n}}J_{n,j}^{(q)} \cap A) - \mu J_{n,j}^{(q)} \mu A \right|
\leq \int \left|\frac{1}{r_{n}+1}\sum_{i=0}^{r_{n}-1} \chi_{A} \circ T^{-i}\right| d\mu
+ \frac{1}{r_{n}+1}
+ \frac{r_{n}}{\tilde{h}_{n}}
\]
and
\[
\sum_{q} \sum_{j=0}^{\tilde{h}_{n}-1}\left|\mu(T^{\tilde{h}_{n}}I_{n,j} \cap A \cap \tilde{C}_{n}^{(q)}) - \mu I_{n,j} \mu(A \cap \tilde{C}_{n}^{(q)}) \right|
\leq \int \left|\frac{1}{r_{n}}\sum_{i=1}^{r_{n}-1} \chi_{A} \circ T^{-i}\right| d\mu
+ \frac{1}{r_{n}}
+ \frac{r_{n}}{h_{n}}
\]
\end{lemma}
\begin{proof}
The first is immediate from the first statement in Lemma \ref{mix0} as the meshes are disjoint.  For the second, observe that $\mu J_{n,j}^{(q)} = \rho \mu I_{n,j}$ where $\rho$ is the size of $\tilde{C}_{n}^{(q)}$ and that, since $A$ is a union of levels, $\mu(A \cap \tilde{C}_{n}^{(q)}) = \rho \mu(A)$.  Then the claim follows from the second statement in Lemma \ref{mix0} and that, as the meshes are disjoint,
\[
\sum_{q} \sum_{j=0}^{\tilde{h}_{n}-1}\mu (T^{\tilde{h}_{n}}I_{n,j}^{[r_{n}]} \cap \tilde{C}_{n}^{(q)})
=  \sum_{j=0}^{\tilde{h}_{n}-1}\mu (T^{\tilde{h}_{n}}I_{n,j}^{[r_{n}]})
= \frac{1}{r_{n}+1}  \sum_{j=0}^{\tilde{h}_{n}-1}\mu (T^{\tilde{h}_{n}}I_{n,j})
\leq \frac{1}{r_{n}+1} \qedhere
\]
\end{proof}

\begin{lemma}\label{mix10}
Let $T$ be a rank-one  transformation and $n$ be such that $B_{n+1}$ is constructed using the staircase operation.  
Let $\tilde{C}_{n} = \bigsqcup_{j=0}^{\tilde{h}_{n}-1}J_{n,j}$ be a mesh for $C_{n}$.
Let $1 \leq k < r_{n}$ and $A$ be a union of levels in $C_{n}$.  Then for $j \geq kr_{n} + \frac{1}{2}k(k-1)$,
\begin{align*}
\left|\mu(T^{k\tilde{h}_{n}}J_{n,j} \cap A) - \mu J_{n,j} \mu A\right|
&\leq \int_{T^{-\frac{1}{2}k(k-1)}J_{n,j}} \left|\frac{1}{r_{n}+1}\sum_{i=0}^{r_{n}-k}\chi_{A} \circ T^{-ik}\right| d\mu \\
&\quad\quad + \left|\sum_{i=r_{n}-k+1}^{r_{n}}\mu(T^{k\tilde{h}_{n}}J_{n,j}^{[i]} \cap A) - \mu J_{n,j}^{[i]}\mu A\right|
\end{align*}
and therefore
\[
\left|\mu(T^{k\tilde{h}_{n}}J_{n,j} \cap A) - \mu J_{n,j} \mu A\right|
\leq \int_{T^{-\frac{1}{2}k(k-1)}J_{n,j}} \left|\frac{1}{r_{n}+1}\sum_{i=0}^{r_{n}-k-1}\chi_{A} \circ T^{-ik}\right| d\mu
+ \frac{k}{r_{n}+1}\mu J_{n,j}
\]
\end{lemma}
\begin{proof}
For $j \geq kr_{n}$ and $0 \leq i \leq r_{n} - k$,
\[
T^{k\tilde{h}_{n}}J_{n,j}^{[i]} = T^{-ik-\frac{1}{2}k(k-1)}J_{n,j}^{[i+k]} = J_{n,j-ik-\frac{1}{2}k(k-1)}^{[i+k]}
\]
so by Lemma \ref{mesh1},
\[
\mu(T^{k\tilde{h}_{n}}J_{n,j}^{[i]} \cap A) = \mu(J_{n,j-ik-\frac{1}{2}k(k-1)}^{[i+k]} \cap A) = \frac{1}{r_{n}+1}\mu(J_{n,j-ik-\frac{1}{2}k(k-1)} \cap A)
= \frac{1}{r_{n}+1}\mu(T^{-ik-\frac{1}{2}k(k-1)}J_{n,j} \cap A)
\]
and therefore
\[
\mu(T^{k\tilde{h}_{n}}J_{n,j} \cap A) = \frac{1}{r_{n}+1}\sum_{i=0}^{r_{n}-k}\mu(T^{-ik-\frac{1}{2}k(k-1)}J_{n,j} \cap A) +
\sum_{i=r_{n}-k+1}^{r_{n}} \mu(T^{k\tilde{h}_{n}}J_{n,j}^{[i]} \cap A)
\]
Then
\begin{align*}
&\left|\mu(T^{k\tilde{h}_{n}}J_{n,j} \cap A) - \mu J_{n,j} \mu A\right| \\
&\quad\quad\leq \left|\int_{T^{-\frac{1}{2}k(k-1)}J_{n,j}} \frac{1}{r_{n}+1}\sum_{i=0}^{r_{n}-k}\chi_{A} \circ T^{-ik}~d\mu\right|
+ \left|\sum_{i=r_{n}-k+1}^{r_{n}}\mu(T^{k\tilde{h}_{n}}J_{n,j}^{[i]} \cap A) - \mu J_{n,j}^{[i]}\mu A\right| \\
&\quad\quad\leq \int_{T^{-\frac{1}{2}k(k-1)}J_{n,j}} \left|\frac{1}{r_{n}+1}\sum_{i=0}^{r_{n}-k}\chi_{A} \circ T^{-ik}\right| d\mu 
+ \left|\sum_{i=r_{n}-k+1}^{r_{n}}\mu(T^{k\tilde{h}_{n}}J_{n,j}^{[i]} \cap A) - \mu J_{n,j}^{[i]}\mu A\right| \qedhere
\end{align*}
\end{proof}

\begin{lemma}\label{mix10inv}
Let $T$ be a rank-one  transformation and $n$ be such that $B_{n+1}$ is constructed using the staircase operation.  
Let $\tilde{C}_{n} = \bigsqcup_{j=0}^{\tilde{h}_{n}-1}J_{n,j}$ be a mesh for $C_{n}$.
Let $1 \leq k < r_{n}$ and $A$ be a union of levels in $C_{n}$.  Then for $j \geq kr_{n} + \frac{1}{2}k(k-1)$,
\begin{align*}
\left|\mu(T^{k\tilde{h}_{n}}I_{n,j} \cap A \cap \tilde{C}_{n}) - \mu I_{n,j} \mu(A \cap \tilde{C}_{n})\right|
&\leq \int_{T^{-\frac{1}{2}k(k-1)}J_{n,j}} \left|\frac{1}{r_{n}+1}\sum_{i=0}^{r_{n}-k}\chi_{A} \circ T^{-ik}\right| d\mu \\
&\quad\quad + \left|\sum_{i=r_{n}-k+1}^{r_{n}}\mu(T^{k\tilde{h}_{n}}I_{n,j}^{[i]} \cap A \cap \tilde{C}_{n}) - \mu I_{n,j}^{[i]}\mu (A \cap \tilde{C}_{n})\right|
\end{align*}
and therefore
\[
\left|\mu(T^{k\tilde{h}_{n}}I_{n,j} \cap A \cap \tilde{C}_{n}) - \mu I_{n,j} \mu (A \cap \tilde{C}_{n})\right|
\leq \int_{T^{-\frac{1}{2}k(k-1)}J_{n,j}} \left|\frac{1}{r_{n}+1}\sum_{i=0}^{r_{n}-k}\chi_{A} \circ T^{-ik}\right| d\mu
+ \frac{k}{r_{n}+1}\mu I_{n,j} \mu \tilde{C}_{n}
\]
\end{lemma}
\begin{proof}
For $j \geq kr_{n}$ and $1 \leq i \leq r_{n} - k$,
\[
T^{k\tilde{h}_{n}}I_{n,j}^{i} = T^{-ik-\frac{1}{2}k(k-1)}I_{n,j}^{[i+k]} = I_{n,j-ik-\frac{1}{2}k(k-1)}^{[i+k]}
\]
so by Lemma \ref{mesh1}, as $I_{n,j-ik-\frac{1}{2}k(k-1)}^{[i+k]} \cap \tilde{C}_{n} = J_{n,j-ik}^{[i+k]}$,
\begin{align*}
\mu(T^{k\tilde{h}_{n}}I_{n,j} \cap A \cap \tilde{C}_{n}) &=
\mu(I_{n,j-ik-\frac{1}{2}k(k-1)}^{[i+k]} \cap A \cap \tilde{C}_{n}) \\
&= \mu(J_{n,j-ik-\frac{1}{2}k(k-1)}^{[i+k]} \cap A) = \frac{1}{r_{n}+1}\mu(J_{n,j-ik} \cap A)
= \mu(T^{-ik-\frac{1}{2}k(k-1)}J_{n,j} \cap A)
\end{align*}
so the claim follows as in the proof of Lemma \ref{mix10} making use of the facts that $\mu J_{n,j} \mu A = \mu I_{n,j} \mu(A \cap \tilde{C}_{n})$ and $\mu J_{n,j} = \mu I_{n,j} \frac{\mu \tilde{C}_{n}}{\mu C_{n}} \leq \mu I_{n,j} \mu \tilde{C}_{n}$.
\end{proof}

\begin{lemma}\label{mix100}
Let $T$ be a rank-one  transformation and $n$ be such that $B_{n+1}$ is constructed using the staircase operation.  Let $\tilde{C}_{n+1} = \bigsqcup_{j=0}^{\tilde{h}_{n+1}-1} J_{n+1,j}$ be a mesh for $C_{n+1}$ (so by Lemma \ref{meshdrop}, $\tilde{C}_{n} = C_{n} \cap \tilde{C}_{n+1} = \bigsqcup_{j=0}^{\tilde{h}_{n}-1}J_{n,j}$ is a mesh for $C_{n}$).  Let $1 \leq k < r_{n}$.  
Then for $j \geq r_{n}(r_{n}-1) + c_{n+1} - c_{n}$ and $r_{n} - k < i \leq r_{n}$,
\[
T^{k\tilde{h}_{n}}J_{n,j}^{[i]}
= T^{\tilde{h}_{n+1}} J_{n+1,j+\frac{1}{2}i(i-1) - \frac{1}{2}r_{n}(r_{n}-1) - c_{n+1} + c_{n}}
\]
\end{lemma}
\begin{proof}
Since $J_{n,j}^{[i]} = J_{n+1,j + i\tilde{h}_{n} + \frac{1}{2} i (i - 1)}$ and $\tilde{h}_{n+1} = (r_{n}+1)\tilde{h}_{n} + \frac{1}{2}r_{n}(r_{n}-1) + c_{n+1} - c_{n}$, this follows from the fact that
\[
i\tilde{h}_{n} + \frac{1}{2} i (i - 1) + k\tilde{h}_{n} = \tilde{h}_{n+1} + (i + k - r_{n} - 1)\tilde{h}_{n} + \frac{1}{2} i (i - 1) - \frac{1}{2}r_{n}(r_{n}-1) - c_{n+1} + c_{n} \qedhere
\]
\end{proof}

\begin{lemma}\label{mix1}
Let $T$ be a rank-one transformation and $n$ be such that $B_{n+1}$ and $B_{n+2}$ are constructed using the staircase operation with $c_{n+1} = c_{n} + r_{n}$.  Let $\tilde{C}_{n+1} = \bigsqcup_{j=0}^{\tilde{h}_{n+1}-1} J_{n+1,j}$ be a mesh for $C_{n+1}$ (so by Lemma \ref{meshdrop}, $\tilde{C}_{n} = C_{n} \cap \tilde{C}_{n+1} = \bigsqcup_{j=0}^{\tilde{h}_{n}-1}J_{n,j}$ is a mesh for $C_{n}$).
Then for $1 \leq k < r_{n}$ and any union of levels $A$ in $C_{n}$,
\begin{align*}
&\sum_{j=0}^{\tilde{h}_{n}-1}\left|\mu(T^{k\tilde{h}_{n}}J_{n,j} \cap A) - \mu J_{n,j} \mu A \right| 
\leq \int_{T^{-\frac{1}{2}k(k-1)}\tilde{C}_{n}} \left|\frac{1}{r_{n}+1}\sum_{i=0}^{r_{n}-k} \chi_{A} \circ T^{-ik}\right| d\mu \\
&\quad\quad\quad\quad + \int_{\tilde{C}_{n+1}}\left|\frac{1}{r_{n+1}+1}\sum_{i=0}^{r_{n+1}-1} \chi_{A}\circ T^{-i}\right| d\mu 
+ \frac{r_{n}^{2}}{\tilde{h}_{n}}\mu \tilde{C}_{n} + \frac{1}{r_{n+1}+1} \mu \tilde{C}_{n+1}
\end{align*}
\end{lemma}
\begin{proof}
For $r_{n} - k < i \leq r_{n}$ and $r_{n}(r_{n}-1) + c_{n+1} - c_{n} \leq j < h_{n}$, by Lemma \ref {mix100} then Lemma \ref{mix00},
\begin{align*}
&\left|\mu(T^{k\tilde{h}_{n}}J_{n,j}^{[i]} \cap A) - \mu J_{n,j}^{[i]} \mu A \right| \\
&\quad\quad= \left|\mu(T^{\tilde{h}_{n+1}}J_{n+1,j + \frac{1}{2}i(i-1) - \frac{1}{2}r_{n}(r_{n}-1) - c_{n+1} + c_{n}} \cap A) - \mu J_{n+1,j + \frac{1}{2}i(i-1) - \frac{1}{2}r_{n}(r_{n}-1)} \mu A \right| \\
&\quad\quad\leq \int_{J_{n+1,j + \frac{1}{2}i(i-1) - \frac{1}{2}r_{n}(r_{n}-1) - c_{n+1} + c_{n}}} \left|\frac{1}{r_{n+1}+1}\sum_{i=0}^{r_{n+1}-2}\chi_{A} \circ T^{-i}\right| d\mu \\
&\quad\quad\quad\quad + \frac{1}{r_{n+1}+1}\mu(J_{n+1,j + \frac{1}{2}i(i-1) - \frac{1}{2}r_{n}(r_{n}-1) - c_{n+1} + c_{n}})
\end{align*}
Then, since the levels in $C_{n+1}$ are disjoint over $j$ and $i$,
\begin{align*}
&\sum_{j=r_{n}(r_{n}-1)+c_{n+1}}^{h_{n}} \left|\sum_{i=r_{n}-k+1}^{r_{n}}\mu(T^{k\tilde{h}_{n}}J_{n,j}^{[i]} \cap A) - \mu J_{n,j}^{[i]}\mu A\right| \\
&\quad\quad \leq \int_{\tilde{C}_{n+1}} \left|\frac{1}{r_{n+1}+1}\sum_{i=0}^{r_{n+1}-1}\chi_{A} \circ T^{-i}\right| d\mu + \frac{1}{r_{n+1}+1}\mu \tilde{C}_{n+1}
\end{align*}
Then by Lemma \ref{mix10},
\begin{align*}
&\sum_{j=r_{n}(r_{n}-1)+c_{n+1}}^{h_{n}} \left|\mu(T^{k\tilde{h}_{n}}J_{n,j} \cap A) - \mu J_{n,j} \mu A\right| \\
&\leq \int_{T^{-\frac{1}{2}k(k-1)}\tilde{C}_{n}} \left|\frac{1}{r_{n}+1}\sum_{i=0}^{r_{n}-k}\chi_{A}\circ T^{-ik}\right| d\mu + \int_{\tilde{C}_{n+1}} \left|\frac{1}{r_{n+1}+1}\sum_{i=0}^{r_{n+1}-1}\chi_{A} \circ T^{-i}\right| d\mu + \frac{1}{r_{n+1}+1}\mu \tilde{C}_{n+1}
\end{align*}
As the bottom $r_{n}(r_{n}-1)+c_{n+1}-c_{n}$ levels in $\tilde{C}_{n}$ have measure at most $\frac{r_{n}(r_{n}-1) + c_{n+1} - c_{n}}{h_{n}}\mu \tilde{C}_{n}$, since $i \leq r_{n} - k$, then $ik + \frac{1}{2}k(k-1) \leq r_{n}k - \frac{1}{2}k^{2} - \frac{1}{2}k < r_{n}(r_{n}-1)$ and $c_{n+1} = c_{n} + r_{n}$, the claim follows.
\end{proof}

\begin{lemma}\label{mix1inv}
Let $T$ be a rank-one transformation and $n$ be such that $B_{n+1}$ and $B_{n+2}$ are constructed using the staircase operation with $c_{n+1} = c_{n} + r_{n}$.  Let $\tilde{C}_{n+1} = \bigsqcup_{j=0}^{\tilde{h}_{n+1}-1} J_{n+1,j}$ be a mesh for $C_{n+1}$ (so by Lemma \ref{meshdrop}, $\tilde{C}_{n} = C_{n} \cap \tilde{C}_{n+1} = \bigsqcup_{j=0}^{\tilde{h}_{n}-1}J_{n,j}$ is a mesh for $C_{n}$).
Then for $1 \leq k < r_{n}$ and any union of levels $A$ in $C_{n}$,
\begin{align*}
&\sum_{j=0}^{\tilde{h}_{n}-1}\left|\mu(T^{k\tilde{h}_{n}}I_{n,j} \cap A \cap \tilde{C}_{n+1}) - \mu I_{n,j} \mu(A \cap \tilde{C}_{n+1}) \right| \\
&\quad\quad \leq \int_{T^{-\frac{1}{2}k(k-1)}\tilde{C}_{n}} \left|\frac{1}{r_{n}+1}\sum_{i=0}^{r_{n}-k} \chi_{A} \circ T^{-ik}\right| d\mu \\
&\quad\quad\quad\quad + \int_{\tilde{C}_{n+1}}\left|\frac{1}{r_{n+1}+1}\sum_{i=0}^{r_{n+1}-1} \chi_{A}\circ T^{-i}\right| d\mu 
+ \frac{r_{n}^{2}}{h_{n}}\mu \tilde{C}_{n} + \frac{1}{r_{n+1}+1} \mu \tilde{C}_{n+1}
\end{align*}
\end{lemma}
\begin{proof}
For $r_{n} - k \leq i < r_{n}$ and $r_{n}(r_{n}-1) + c_{n+1} - c_{n} \leq j < h_{n}$, by Lemma \ref {mix100} then Lemma \ref{mix00inv}, and using that $\mu I_{n+1,j} \mu(A \cap \tilde{C}_{n+1}) = \mu J_{n+1,j} \mu A$,
\begin{align*}
&\left|\mu(T^{k\tilde{h}_{n}}I_{n,j}^{[i]} \cap A \cap \tilde{C}_{n+1}) - \mu I_{n,j}^{[i]} \mu (A \cap \tilde{C}_{n+1}) \right| \\
&\quad\quad = \left|\mu(T^{\tilde{h}_{n+1}}I_{n+1,j + \frac{1}{2}i(i-1) - \frac{1}{2}r_{n}(r_{n}-1) - c_{n+1} + c_{n}} \cap A \cap \tilde{C}_{n+1}) - \mu J_{n+1,j + \frac{1}{2}i(i-1) - \frac{1}{2}r_{n}(r_{n}-1) - c_{n+1} + c_{n}} \mu A \right| \\
&\quad\quad \leq \int_{J_{n+1,j + \frac{1}{2}i(i-1) - \frac{1}{2}r_{n}(r_{n}-1) - c_{n+1} + c_{n}}} \left|\frac{1}{r_{n+1}+1}\sum_{i=0}^{r_{n+1}-2}\chi_{A} \circ T^{-i}\right| d\mu  \\
&\quad\quad\quad\quad + \frac{1}{r_{n+1}+1}\mu(J_{n+1,j + \frac{1}{2}i(i-1) - \frac{1}{2}r_{n}(r_{n}-1) - c_{n+1} + c_{n}})
\end{align*}
and then the claim follows as in the proof of Lemma \ref{mix1}.
\end{proof}

\begin{lemma}\label{mix1all}
Let $T$ be a rank-one transformation and $n$ be such that $B_{n+1}$ and $B_{n+2}$ are constructed using the staircase operation with $c_{n+1} = c_{n} + r_{n}$.  Let $\tilde{C}_{n+1}^{(q)}$ be a collection of disjoint meshes of $C_{n+1}$ and $A$ be a union of levels in $C_{n}$.  Then for $1 \leq k < r_{n}$,
\begin{align*}
\sum_{q} \sum_{j=0}^{\tilde{h}_{n}-1}&\left|\mu(T^{k\tilde{h}_{n}}J_{n,j}^{(q)} \cap A) - \mu J_{n,j}^{(q)} \mu A \right| \\
&\leq \int \left|\frac{1}{r_{n}+1}\sum_{i=0}^{r_{n}-k} \chi_{A} \circ T^{-ik}\right| d\mu
+ \int \left|\frac{1}{r_{n+1}+1}\sum_{i=0}^{r_{n+1}} \chi_{A}\circ T^{-i}\right| d\mu 
+ \frac{r_{n}^{2}}{\tilde{h}_{n}} + \frac{1}{r_{n+1}+1} 
\end{align*}
and
\begin{align*}
\sum_{q} \sum_{j=0}^{\tilde{h}_{n}-1}&\left|\mu(T^{k\tilde{h}_{n}}I_{n,j} \cap A \cap \tilde{C}_{n+1}) - \mu I_{n,j} \mu(A \cap \tilde{C}_{n+1}) \right| \\
&\leq \int \left|\frac{1}{r_{n}+1}\sum_{i=0}^{r_{n}-k} \chi_{A} \circ T^{-ik}\right| d\mu
+ \int \left|\frac{1}{r_{n+1}+1}\sum_{i=0}^{r_{n+1}-1} \chi_{A}\circ T^{-i}\right| d\mu 
+ \frac{r_{n}^{2}}{\tilde{h}_{n}} + \frac{1}{r_{n+1}+1} 
\end{align*}
\end{lemma}
\begin{proof}
This is immediate from Lemmas \ref{mix1} and \ref{mix1inv} using the disjointness.
\end{proof}

\begin{lemma}\label{mix2}
Let $T$ be a rank-one transformation and $n$ be such that $B_{n+1}$ is constructed using the staircase operation and $B_{n+2}$ through $B_{n+z}$ using Chacon-like operations with the same $M_{n+1}$ and $c_{n+1} < c_{n} + r_{n} + M_{n+1}$ and $B_{n+z+1}$ be constructed using the staircase operation.  Let $\tilde{C}_{n}^{(q)} = \bigsqcup_{j=0}^{\tilde{h}_{n}-1}J_{n,j}^{(q)}$ be disjoint meshes for $C_{n}$.
Then for $1 \leq k < r_{n}$ and any union of levels $A$ in $C_{n}$,
\begin{align*}
\sum_{q} \sum_{j=0}^{\tilde{h}_{n}-1}\left|\mu(T^{k\tilde{h}_{n}}J_{n,j}^{(q)} \cap A) - \mu J_{n,j}^{(q)} \mu A \right| 
&\leq \int \left|\frac{1}{r_{n}+1}\sum_{i=0}^{r_{n}-k} \chi_{A} \circ T^{-ik}\right| d\mu
+ \int \left|\frac{1}{r_{n}+1}\sum_{i=0}^{k} \chi_{A} \circ T^{i(r_{n}-k)}\right| d\mu \\
&+ \int \left|\frac{1}{r_{n+z}+1}\sum_{i=0}^{r_{n+z}-1}\chi_{A}\circ T^{-i}\right| d\mu + \frac{1}{r_{n+z}} + \frac{r_{n}^{2} + c_{n} + 2M_{n+1}}{h_{n}}
\end{align*}
\end{lemma}
\begin{proof}
Let $\tilde{C}_{n} = \bigsqcup_{j=0}^{\tilde{h}_{n-1}} J_{n,j}$ be one of the meshes.
Let $S^{\star}$ be the spacer levels added above $C_{n+z}$.
Let $Q = \{ x \in I_{n+1,0} : \exists~0 < y \leq \tilde{h}_{n+1}~\text{s.t.}~T^{y}x \in S^{\star} \}$.
Let $L_{n+1} = \{ x \in I_{n+1,0} \setminus Q : T^{\tilde{h}_{n+1}}x \in I_{n+1,0} \}$ and $L_{n+1}^{\star} = \{ x \in I_{n+1,0} \setminus Q : T^{\tilde{h}_{n+1}+M_{n+1}}x \in I_{n+1,0} \}$.  Since $B_{n+2}$ through $B_{n+z}$ are constructed using Chacon-like operations all with the same $M_{n+1}$ and $c_{n+1}$, $I_{n+1,0} = L_{n+1} \sqcup L_{n+1}^{\star} \sqcup Q$.

Let $L_{n} = \bigsqcup_{i=1}^{r_{n}} T^{i\tilde{h}_{n} + \frac{1}{2}i(i-1)}L_{n+1}$ and $L_{n}^{\star} = \bigsqcup_{i=1}^{r_{n}} T^{(i\tilde{h}_{n} + \frac{1}{2}i(i-1)}L_{n+1}^{\star}$.  For $x \in L_{n+1}$, since $x \notin Q$ and $T^{\tilde{h}_{n+1}}x \in I_{n+1,0}$, it follows that if $x \in T^{i\tilde{h}_{n} + \frac{1}{2}i(i-1)}L_{n+1}$ then $T^{\tilde{h}_{n+1}}x \in T^{i\tilde{h}_{n} + \frac{1}{2}i(i-1)}L_{n+1}$.  Let $K_{n,0} = T^{\tilde{h}_{n+1}}(J_{n,0} \cap L)$.  Then $K_{n,0} \subseteq I_{n,0}$ and $K_{n,0}^{[i]} = T^{i\tilde{h}_{n}+\frac{1}{2}i(i-1)}K_{n,0}^{[0]}$ so $K_{n,0}$ is a mesh prototype for $C_{n}$.  Likewise, $K_{n,0}^{\star} = T^{\tilde{h}_{n+1}}(J_{n,0} \cap L^{\star})$ is a mesh prototype for $C_{n}$.

For $r_{n} - k < i \leq r_{n}$ and $0 \leq j < \tilde{h}_{n}$, since $\tilde{h}_{n+1} = (r_{n}+1)\tilde{h}_{n} + \frac{1}{2}r_{n}(r_{n}-1) + c_{n+1} - c_{n}$,
\begin{align*}
T^{k\tilde{h}_{n}}J_{n,j}^{[i]}
&= T^{k\tilde{h}_{n}}T^{(i\tilde{h}_{n}+\frac{1}{2}i(i-1)} J_{n,j}^{[0]}
= T^{\tilde{h}_{n+1}}T^{(i+k - r_{n}-1)\tilde{h}_{n} + \frac{1}{2}i(i-1) - \frac{1}{2}r_{n}(r_{n}-1) - c_{n+1} + c_{n}} J_{n,j}^{[0]} \\
&= T^{\tilde{h}_{n+1}}T^{\frac{1}{2}i(i-1) - \frac{1}{2}r_{n}(r_{n}-1) - c_{n+1} + c_{n} - \frac{1}{2}(i+k-r_{n})(i+k-r_{n}-1)}J_{n,j}^{[i+k-r_{n}]} \\
&= T^{\tilde{h}_{n+1}}T^{i(r_{n}-k) - x_{n,k}}J_{n,j}^{[i+k-r_{n}]}
\end{align*}
where $x_{n,k} =  r_{n}(r_{n}-1) - \frac{1}{2}k(r_{n}-k+1) + c_{n+1} - c_{n}$.  Since $c_{n+1} < c_{n} + r_{n} + M_{n+1}$,
$x_{n,k} < r_{n}^{2} + M_{n+1}$

Then for $r_{n}^{2} + M_{n+1} \leq j < \tilde{h}_{n}$,
\[
T^{k\tilde{h}_{n}}J_{n,j}^{[i]}
= T^{\tilde{h}_{n+1}}J_{n,j+i(r_{n}-k)-x_{n,k}}^{[i+k-r_{n}]}
\]
Now observe that
\begin{align*}
T^{\tilde{h}_{n+1}}J_{n,j} = T^{\tilde{h}_{n+1}}T^{j}J_{n,0}
&= T^{j}\left(T^{\tilde{h}_{n+1}}\left((J_{n,0} \cap L_{n}) \sqcup (J_{n,0} \cap L_{n}^{\star}) \sqcup (J_{n,0} \cap Q)\right)\right) \\
&= T^{j}\left(K_{n,0} \sqcup T^{-M_{n+1}}K_{n,0}^{\star} \sqcup T^{\tilde{h}_{n+1}}(J_{n,0} \cap Q)\right)
\end{align*}
Since $J_{n,0}$, $K_{n,0}$ and $K_{n,0}^{\star}$ are all meshes, then for $j \geq r_{n}^{2} + 2M_{n+1}$,
\[
T^{\tilde{h}_{n+1}}J_{n,j+i(r_{n}-k)-x_{n,k}}^{[i+k-r_{n}]} = K_{n,j+i(r_{n}-k)-x_{n,k}}^{[i+k-r_{n}]} \sqcup  (K_{n,j+i(r_{n}-k)-x_{n,k}-M_{n+1}}^{\star})^{[i+k-r_{n}]} \sqcup T^{\tilde{h}_{n+1}}Z_{i,j}
\]
where $Z_{i,j}$ is a collection of levels in $C_{n+z}$ such that $T^{\tilde{h}_{n+1}}$ pushes each level through the top of $C_{n+z}$.

Since $K_{n,0}$ and $K_{n,0}^{\star}$ are mesh prototypes for $C_{n}$ and $A$ is a union of levels in $C_{n}$, then by Lemma \ref{mesh1},
\[
\mu(K_{n,j+i(r_{n}-k)-x_{n,k}}^{[i+k-r_{n}]} \cap A)
= \frac{1}{r_{n}+1}\mu(K_{n,j+i(r_{n}-k)-x_{n,k}} \cap A)
= \frac{1}{r_{n}+1}\mu(T^{i(r_{n}-k) - x_{n,k}}K_{n,j} \cap A)
\]
and
\begin{align*}
\mu((K_{n,j+i(r_{n}-k)-x_{n,k}-M_{n+1}}^{\star})^{[i+k-r_{n}]} \cap A)
&= \frac{1}{r_{n}+1}\mu(K_{n,j+i(r_{n}-k)-x_{n,k}-M_{n+1}} \cap A) \\
&= \frac{1}{r_{n}+1}\mu(T^{i(r_{n}-k) - x_{n,k} - M_{n+1}}K_{n,j} \cap A)
\end{align*}
Therefore
\begin{align*}
&\left|\sum_{i=r_{n}-k+1}^{r_{n}} \mu(T^{k\tilde{h}_{n}}J_{n,j}^{[i]} \cap A) - \mu J_{n,j}^{[i]}\mu A\right|
\leq \left|\frac{1}{r_{n}+1}\sum_{i=r_{n}-k+1}^{r_{n}}\mu(T^{i(r_{n}-k) - x_{n,k}}K_{n,j} \cap A) - \mu K_{n,j}\mu A\right| \\
&\quad\quad + \left|\frac{1}{r_{n}+1}\sum_{i=r_{n}-k+1}^{r_{n}}\mu(T^{i(r_{n}-k) - x_{n,k} - M_{n+1}}K_{n,j} \cap A) - \mu K_{n,j} \mu A \right| \\
&\quad\quad +  \sum_{\text{levels in $C_{n+z}$}} \mu(T^{\tilde{h}_{n+z}}I_{n+z,y} \cap A)
\end{align*}
Using that $T$ is measure-preserving,
\begin{align*}
&\left|\frac{1}{r_{n}+1}\sum_{i=r_{n}-k+1}^{r_{n}}\mu(T^{i(r_{n}-k) - x_{n,k}}K_{n,j} \cap A) - \mu K_{n,j}\mu A\right| \\
&\quad\quad = \left|\frac{1}{r_{n}+1}\sum_{i=1}^{k}\mu(T^{i(r_{n}-k) + (r_{n}-k)^{2} - x_{n,k}} K_{n,j} \cap A) - \mu K_{n,j}\mu A\right| \\
&\quad\quad\leq \int_{T^{(r_{n}-k)^{2}-x_{n,k}} K_{n,j}} \left|\frac{1}{r_{n}+1}\sum_{i=0}^{k}\chi_{A} \circ T^{i(r_{n}-k)}\right|d\mu
\end{align*}
and similarly
\begin{align*}
&\left|\frac{1}{r_{n}+1}\sum_{i=r_{n}-k+1}^{r_{n}}\mu(T^{i(r_{n}-k) - x_{n,k} - M_{n+1}}K_{n,j} \cap A) - \mu K_{n,j}^{\star} \mu A \right| \\
&\quad\quad \leq \int_{T^{(r_{n}-k)^{2}-x_{n,k}-M_{n+1}} K_{n,j}^{\star}} \left|\frac{1}{r_{n}+1}\sum_{i=1}^{k}\chi_{A} \circ T^{i(r_{n}-k)}\right|d\mu
\end{align*}
Since $K_{n,0}$ and $T^{-M_{n}}K_{n,0}^{\star}$ are disjoint and contained in $T^{\tilde{h}_{n+1}}J_{n,0}$, then
\begin{align*}
&\int_{T^{(r_{n}-k)^{2}-x_{n,k} K_{n,j}}} \left|\frac{1}{r_{n}+1}\sum_{i=0}^{k}\chi_{A} \circ T^{i(r_{n}-k)}\right|d\mu + \int_{T^{(r_{n}-k)^{2}-x_{n,k}-M_{n+1} K_{n,j}^{\star}}} \left|\frac{1}{r_{n}+1}\sum_{i=0}^{k}\chi_{A} \circ T^{i(r_{n}-k)}\right|d\mu \\
&\quad\quad \leq \int_{T^{(r_{n}-k)^{2}-x_{n,k}+\tilde{h}_{n+1}} J_{n,j}} \left|\frac{1}{r_{n}+1}\sum_{i=0}^{k}\chi_{A} \circ T^{i(r_{n}-k)}\right|d\mu 
\end{align*}
Therefore if $\tilde{C}_{n}^{(q)}$ are disjoint meshes of $C_{n}$ then by Lemma \ref{mix0all} applied to the levels of $C_{n+z}$ pushed through the top of $C_{n+z}$,
\begin{align*}
\sum_{q} \left|\sum_{i=r_{n}-k+1}^{r_{n}} \mu(T^{k\tilde{h}_{n}} (J_{n,j}^{(q)})^{[i]} \cap A) - \mu (J_{n,j}^{(q)})^{[i]}\mu A\right|
&\leq \int \left|\frac{1}{r_{n}+1}\sum_{i=0}^{k}\chi_{A} \circ T^{i(r_{n}-k)}\right|d\mu 
+ \frac{r_{n}^{2} + 2M_{n+1}}{h_{n}} \\
&\quad\quad + \int \left|\frac{1}{r_{n+z}+1}\sum_{i=0}^{r_{n+z}-1}\chi_{A} \circ T^{-i}\right|d\mu 
+ \frac{1}{r_{n+z}}
\end{align*}
so the claim follows from Lemma \ref{mix10inv}.
\end{proof}

\begin{lemma}\label{mix2inv}
Let $T$ be a rank-one transformation and $n$ be such that $B_{n+1}$ is constructed using the staircase operation and $B_{n+2}$ through $B_{n+z}$ using Chacon-like operations with the same $M_{n+1}$ and $c_{n+1} < c_{n} + r_{n} + M_{n+1}$ and $B_{n+z+1}$ be constructed using the staircase operation.  Let $\tilde{C}_{n}^{(q)} = \bigsqcup_{j=0}^{\tilde{h}_{n}-1}J_{n,j}^{(q)}$ be disjoint meshes for $C_{n}$.
Then for $1 \leq k < r_{n}$ and any union of levels $A$ in $C_{n}$,
\begin{align*}
\sum_{q} \sum_{j=0}^{\tilde{h}_{n}-1}&\left|\mu(T^{k\tilde{h}_{n}}I_{n,j} \cap A \cap \tilde{C}_{n}^{(q)}) - \mu I_{n,j} \mu (A \cap \tilde{C}_{n}^{(q)}) \right| \\
&\leq \int \left|\frac{1}{r_{n}+1}\sum_{i=0}^{r_{n}-k} \chi_{A} \circ T^{-ik}\right| d\mu
+ \int \left|\frac{1}{r_{n}+1}\sum_{i=0}^{k} \chi_{A} \circ T^{i(r_{n}-k)}\right| d\mu \\
&\quad\quad + \int \left|\frac{1}{r_{n+z}}\sum_{i=1}^{r_{n+z}-1}\chi_{A}\circ T^{-i}\right| d\mu + \frac{1}{r_{n+z}} + \frac{r_{n}^{2} + 2M_{n+1}}{\tilde{h}_{n}}
\end{align*}
\end{lemma}
\begin{proof}
Let $\tilde{C}_{n} = \bigsqcup_{j=0}^{\tilde{h}_{n}-1} J_{n,j}$ be one of the meshes.
Let $L_{n,0}$ and $L_{n,0}^{\star}$ and $Q$ be as in the proof of Lemma \ref{mix2}.  Then $L_{n,0}$ and $L_{n,0}^{\star}$ are mesh prototypes.  As in the proof of Lemma \ref{mix2}, for $r_{n} - k < i \leq r_{n}$ and $r_{n}^{2} + 2M_{n+1} \leq j < \tilde{h}_{n+1}$,
\[
T^{k\tilde{h}_{n}}I_{n,j}^{[i]}
= T^{\tilde{h}_{n+1}}I_{n,j+i(r_{n}-k)-x_{n,k}}^{[i+k-r_{n}]}
= L_{n,j+i(r_{n}-k)-x_{n,k}}^{[i+k-r_{n}]} \sqcup (L_{n,j+i(r_{n}-k)-x_{n,k}-M_{n+1}}^{\star})^{[i+k-r_{n}]}
\sqcup T^{\tilde{h}_{n+1}}Z_{i,j}
\]
where $Z_{i,j}$ is a collection of levels in $C_{n+z}$ such that $T^{\tilde{h}_{n+1}}$ pushes each level through the top of $C_{n+z}$.  Since $L_{n,0}$ is a mesh prototype and $\tilde{C}_{n}$ is a mesh, $\tilde{C}_{n} \cap \bigsqcup L_{n,j}$ is a mesh (and likewise with $L_{n,0}^{\star}$) so by Lemma \ref{mesh1},
\begin{align*}
\mu(L_{n,j+i(r_{n}-k)-x_{n,k}}^{[i+k-r_{n}]} \cap A \cap \tilde{C}_{n})
&= \frac{1}{r_{n}+1}\mu(L_{n,j+i(r_{n}-k)-x_{n,k}} \cap A \cap \tilde{C}_{n}) \\
&= \frac{1}{r_{n}+1}\mu(T^{i(r_{n}-k) - x_{n,k}}(L_{n,j} \cap J_{n,j}) \cap A)
\end{align*}
and similarly
\[
\mu((L_{n,j+i(r_{n}-k)-x_{n,k}-M_{n+1}}^{\star})^{[i+k-r_{n}]} \cap A \cap \tilde{C}_{n})
=  \frac{1}{r_{n}+1}\mu(T^{i(r_{n}-k) - x_{n,k} - M_{n+1}}(L_{n,j}^{\star} \cap J_{n,j}) \cap A)
\]
The claim now follows similar to the proof of Lemma \ref{mix2}.
\end{proof}

\begin{lemma}\label{heightmix}
Let $T$ be a rank-one transformation such that the staircase operation is used infinitely often and $r_{n} \to \infty$ along the sequence of $n$ for which the staircase operation is used.  Then $\{ \tilde{h}_{n} \}$ is a mixing sequence along those $n$ for which the staircase operation is used.
\end{lemma}
\begin{proof}
Let $A$ and $D$ be unions of levels in some column $C_{n_{0}}$.  Then for $n \geq n_{0}$ such that the staircase operation is used, as $A$ and $D$ can be written as unions of levels in $C_{n}$, by Lemma \ref{mix0all} applied to the single mesh $\bigsqcup I_{n,j}$,
\begin{align*}
\left|\mu(T^{\tilde{h}_{n}}D \cap A) - \mu D \mu A\right|
&\leq \sum_{j=0}^{\tilde{h}_{n}-1} \left|\mu(T^{\tilde{h}_{n}}I_{n,j} \cap A) - \mu I_{n,j} \mu A \right| \\
&\leq \int\left|\frac{1}{r_{n}+1}\sum_{i=0}^{r_{n}-1}\chi_{A}\circ T^{-i} - \mu A\right|d\mu + \frac{1}{r_{n}+1} + \frac{r_{n}}{\tilde{h}_{n}}
\end{align*}
which tends to zero as $T$ is ergodic and $r_{n} \to \infty$ (and $\frac{r_{n}}{\tilde{h}_{n}} \to 0$ since $T$ is on a finite measure space).
\end{proof}

\begin{lemma}\label{kheightmix}
Let $T$ be a totally ergodic rank-one transformation such that the staircase operation is used infinitely often and $r_{n} \to \infty$ along the sequence of $n$ for which the staircase operation is used.  Then for each fixed $k$, the sequence $\{ k\tilde{h}_{n} \}$ is a mixing sequence along those $n$ for which the staircase operation is used.
\end{lemma}
\begin{proof}
Let $k \geq 1$ be fixed.  Let $A$ and $D$ be unions of levels in some column $C_{n_{0}}$.  Then for $n \geq n_{0}$ such that the staircase operation is used, as $A$ and $D$ can be written as unions of levels in $C_{n}$, by Lemma \ref{mix10} (the second statement) applied to the mesh $\bigsqcup I_{n,j}$,
\begin{align*}
&\left|\mu(T^{k\tilde{h}_{n}}D \cap A) - \mu D \mu A\right| \\
&\quad\quad \leq \sum_{j=0}^{\tilde{h}_{n}-1} \left|\mu(T^{k\tilde{h}_{n}}I_{n,j} \cap A) - \mu I_{n,j} \mu A \right| \\
&\quad\quad \leq \sum_{j=kr_{n} + \frac{1}{2}k(k-1)}^{\tilde{h}_{n}-1} \left( \int_{T^{-\frac{1}{2}k(k-1)}I_{n,j}} \left|\frac{1}{r_{n}+1}\sum_{i=0}^{r_{n}-k}\chi_{A}\circ T^{-ik}\right|d\mu + \frac{k}{r_{n}}\mu I_{n,j} \right) + \frac{kr_{n} + \frac{1}{2}k(k-1)}{\tilde{h}_{n}} \\
&\quad\quad \leq \int \left|\frac{1}{r_{n}+1}\sum_{i=0}^{r_{n}-k}\chi_{A}\circ T^{-ik}\right|d\mu + \frac{k}{r_{n}+1} + \frac{kr_{n} + \frac{1}{2}k(k-1)}{\tilde{h}_{n}}
\end{align*}
which tends to zero as $k$ is fixed, $r_{n} \to \infty$, $\frac{r_{n}}{h_{n}} \to 0$ and $T^{k}$ is ergodic.
\end{proof}

\begin{lemma}[\cite{adams1998smorodinsky}; see Lemma A.16 \cite{CPR}]\label{block}
Let $T$ be a transformation and $n,p,k,L$ be nonzero integers with $n,L > 0$ and $\chi \in L^{1}$.  Then
\[
\int\left|\frac{1}{n}\sum_{i=1}^{n}\chi \circ T^{-ik}\right|d\mu
\leq \int\left|\frac{1}{L}\sum_{i=0}^{L}\chi\circ T^{-ipk}\right|d\mu + \frac{|p|L}{n}
\]
\end{lemma}

\begin{lemma}[\cite{adams1998smorodinsky}; see Lemma A.15 \cite{CPR}]\label{BH}
Let $T$ be a transformation and $J$ a measurable set.
Let $\epsilon > 0$ and $L$ such that $\frac{1}{L} < \epsilon$.  Let $p$ such that for all $1 \leq i < L$,
\[
\left|\mu(T^{ip}J \cap J) - \mu J \mu J \right| < \epsilon
\]
Then
\[
\int \left|\frac{1}{L}\sum_{i=0}^{L-1}\chi_{J} \circ T^{-ip}\right| d\mu < \sqrt{2\epsilon}
\]
\end{lemma}

\begin{definition}
Let $T$ be a rank-one transformation and $A$ be a union of levels in some column $C_{n_{0}}$.  Let $\epsilon > 0$ and $n \in \mathbb{N}$.  Then $A$ is \textbf{uniformly mixed within $\epsilon$ at stage $n$} if for any collection of disjoint meshes $\tilde{C}_{n}^{(q)}$ of $C_{n}$ and all $(1-\frac{1}{2}\epsilon)\tilde{h}_{n} \leq t < (1 - \frac{1}{2}\epsilon)\tilde{h}_{n+1}$,
\[
\sum_{q} \sum_{j=0}^{\tilde{h}_{n}-1}\left|\mu(T^{t}J_{n,j}^{(q)} \cap A) - \mu J_{n,j}^{(q)} \mu A \right| < \epsilon
\]
and
\[
\sum_{q} \sum_{j=0}^{\tilde{h}_{n}-1}\left|\mu(T^{t}I_{n,j} \cap A \cap \tilde{C}_{n}^{(q)}) - \mu I_{n,j} \mu (A \cap \tilde{C}_{n}^{(q)}) \right| < \epsilon
\]
\end{definition}

\begin{lemma}\label{ineq}
Let $T$ be a rank-one transformation and $A$ be a union of levels in some column $C_{n_{0}}$.  Let $\epsilon > 0$ and $n \in \mathbb{N}$ and let $A$ be uniformly mixed within $\epsilon$ at stage $n$.  Then for any collection of disjoint meshes $\tilde{C}_{n}^{(q)}$ of $C_{n}$ and any union of levels $D$ in $C_{n}$ and all $(1-\frac{1}{2}\epsilon)\tilde{h}_{n} \leq t < (1 - \frac{1}{2}\epsilon)\tilde{h}_{n+1}$,
\[
\sum_{q} \left|\mu(T^{t}(D \cap \tilde{C}_{n}^{(q)}) \cap A) - \mu(D \cap \tilde{C}_{n}^{(q)}) \mu A\right| < \epsilon
\]
\end{lemma}
\begin{proof}
Since $D$ is a union of levels in $C_{n}$, there exists $\mathcal{D} \subseteq \{ 0, 1, \ldots, \tilde{h}_{n}-1 \}$ such that $D = \bigsqcup_{j\in\mathcal{D}}I_{n,j}$.  
Then
\begin{align*}
\sum_{q} \left|\mu(T^{t}(D \cap \tilde{C}_{n}^{(q)}) \cap A) - \mu(D \cap \tilde{C}_{n}^{(q)}) \mu A\right| 
&\leq \sum_{q} \sum_{j \in \mathcal{D}} \left|\mu(T^{t}(I_{n,j} \cap \tilde{C}_{n}^{(q)}) \cap A) - \mu(I_{n,j} \cap \tilde{C}_{n}^{(q)}) \mu A \right| \\
&\leq  \sum_{q} \sum_{j=0}^{\tilde{h}_{n}-1} \left|\mu(T^{t}J_{n,j}^{(q)} \cap A) - \mu J_{n,j}^{(q)} \mu A \right| < \epsilon \qedhere
\end{align*}
\end{proof}

\begin{proposition}\label{mix4}
Let $T$ be an elevated staircase-Chacon-like hybrid transformation on a finite measure space such that the staircase operation is used infinitely often and that along the sequence of $n$ for which $B_{n+1}$ is constructed using the staircase operation,
\begin{itemize}
\item $r_{n} \to \infty$
\item $\frac{r_{n}^{2}}{h_{n}} \to 0$
\item $r_{n} = \tilde{h}_{\ell}$ for $\ell < n$ such that $B_{\ell+1}$ was constructed using a staircase operation
\item $c_{n+1} = c_{n} + r_{n}$
\end{itemize}
and along the sequence of $n$ for which $B_{n+1}$ is constructed using a Chacon-like operation,
\begin{itemize}
\item $c_{n} < c_{n-1} + r_{n-1} + M_{n}$
\item for all consecutive Chacon-like operations prior to the next staircase operation, the same $c_{n}$ and $M_{n}$ are used.
\end{itemize}
and such that $\frac{M_{n}}{h_{n-1}} \to 0$ along the sequence of $n$ for which $B_{n+1}$ is constructed using a Chacon-like operation and $B_{n}$ was constructed using a staircase operation.

Let $\epsilon_{0} > 0$ and $A$ be a union of levels in some column $C_{n_{0}}$.  Then there exists $n_{1} > n_{0}$ such that for all $n \geq n_{1}$ for which $B_{n+1}$ is constructed using the staircase operation, $A$ is uniformly mixed within $\epsilon_{0}$ at stage $n$.
\end{proposition}
\begin{proof}
By Lemma \ref{heightmix}, the sequence $\{ \tilde{h}_{n} \}$ restricted to those $n$ for which the staircase operation is used forms a mixing sequence so $T$ is weakly mixing.  Then $T$ is totally ergodic so by Lemma \ref{kheightmix}, for each fixed $k$, again along the staircase heights, $\{ k \tilde{h}_{n} \}$ is mixing.

Choose $\epsilon$ such that $4\sqrt{\epsilon} + 3\sqrt{2\epsilon} + 5\epsilon < \frac{1}{4}\epsilon_{0}$.
Choose $L \geq 1$ such that $\frac{1}{L} < \epsilon$.  Take $n$ such that the staircase operation used to construct $B_{n+1}$ to be large enough that
$\int\left|\frac{1}{r_{n}}\sum_{i=0}^{r_{n}-k}\chi\circ T^{-ik}\right| d\mu < \epsilon$ for all $1 \leq k \leq L^{2}$ and that
$\frac{1}{r_{n}} < \epsilon$ and $\frac{r_{n}^{2} +2M_{n+1}}{h_{n}} < \epsilon$ (or just $\frac{r_{n}^{2}}{h_{n}}$ if the following stage is also a staircase operation)

Since $\{ k\tilde{h}_{n} \}$ is a mixing sequence along $n$ for which $B_{n+1}$ is constructed using the staircase operation, by Lemma \ref{ts}, we may also take $n$ large enough that for all $\tilde{h}_{\ell} \leq t \leq 4L\tilde{h}_{\ell}$,
\[
\left|\mu(T^{t}A \cap A) - \mu A \mu A \right| < \epsilon
\]
since $\ell \to \infty$ as $n \to \infty$ because $r_{n} \to \infty$.

By Lemma \ref{mix1all}, then for $1 \leq k \leq L^{2}$,
\[
\sum_{q} \sum_{j=0}^{\tilde{h}_{n}-1} \left|\mu(T^{k\tilde{h}_{n}}J_{n,j}^{(q)} \cap A) - \mu J_{n,j}^{(q)} \mu A \right| 
< 4\epsilon < \frac{1}{4}\epsilon_{0}
\]
and
\[
\sum_{q} \sum_{j=0}^{\tilde{h}_{n}-1} \left|\mu(T^{k\tilde{h}_{n}}I_{n,j} \cap A \cap \tilde{C}_{n}^{(q)}) - \mu I_{n,j} \mu (A \cap \tilde{C}_{n}^{(q)}) \right| 
< 4\epsilon < \frac{1}{4}\epsilon_{0}
\]

Let $k$ such that $L^{2} \leq k < r_{n} = h_{\ell}$.  Choose $p$ to be the minimal integer such that $pk \geq \tilde{h}_{\ell}$.  Then $pk < 2\tilde{h}_{\ell}$.  Write $\chi = \chi_{A}$.
By Lemma \ref{block},
\[
\int \left|\frac{1}{r_{n}+1}\sum_{i=0}^{r_{n}-1} \chi \circ T^{-ik}\right| d\mu 
\leq \int \left|\frac{1}{L}\sum_{i=0}^{L-1} \chi \circ T^{-ipk}\right| d\mu + \frac{pL}{r_{n}+1}
\]
Since $pk < 2\tilde{h}_{\ell}$ and $r_{n} = \tilde{h}_{\ell}$ and $k \geq L^{2}$,
\[
\frac{pL}{r_{n}+1} <
\frac{pL}{r_{n}} = \frac{pkL}{\tilde{h}_{\ell}k} < \frac{2\tilde{h}_{\ell}L}{\tilde{h}_{\ell}k} = \frac{2L}{k} \leq \frac{2L}{L^{2}} = \frac{2}{L} < 2\epsilon
\]
Now $\tilde{h}_{\ell} \leq ipk < 2L\tilde{h}_{\ell}$ for all $1 \leq i \leq L$ so by Lemma \ref{BH},
\[
\int \left|\frac{1}{L}\sum_{i=0}^{L-1} \chi \circ T^{-ipk}\right| d\mu < \sqrt{2\epsilon} + 2\epsilon
\]
Therefore if $B_{n+2}$ is constructed using the staircase operation, then by Lemma \ref{mix1all},
\[
\sum_{q} \sum_{j=0}^{\tilde{h}_{n}-1} \left|\mu(T^{k\tilde{h}_{n}}J_{n,j}^{(q)} \cap A) - \mu J_{n,j}^{(q)} \mu A \right| 
< \sqrt{2\epsilon} + 2\epsilon + \epsilon + \epsilon + \epsilon < \frac{1}{4}\epsilon_{0}
\]
and
\[
\sum_{q} \sum_{j=0}^{\tilde{h}_{n}-1} \left|\mu(T^{k\tilde{h}_{n}}I_{n,j} \cap A \cap \tilde{C}_{n}^{(q)}) - \mu I_{n,j} \mu (A \cap \tilde{C}_{n}^{(q)}) \right| 
< \frac{1}{4}\epsilon_{0}
\]
Now consider when $B_{n+2}$ through $B_{n+z}$ are constructed using the same Chacon-like operation and $B_{n+z+1}$ using the staircase operation.

If $k < \sqrt{\epsilon} r_{n}$ then $\int \left|\frac{1}{r_{n}}\sum_{i=1}^{k}\chi \circ T^{i(r_{n}-k)}\right| d\mu < \sqrt{\epsilon}$ so we may assume $k \geq \sqrt{\epsilon} r_{n}$.  Likewise, if $r_{n} - k < \frac{1}{2}L^{2}$ then since we have taken $n$ large enough already, $\int \left|\frac{1}{r_{n}}\sum_{i=1}^{k}\chi \circ T^{i(r_{n}-k)}\right| d\mu < \epsilon$.

Take $p$ minimal such that $p(r_{n}-k) \geq \tilde{h}_{\ell}$.  Then $p(r_{n}-k) < 2\tilde{h}_{\ell}$ and
\[
\frac{pL}{k} < \frac{pL}{\sqrt{\epsilon}r_{n}} = \frac{p(r_{n}-k)L}{\sqrt{\epsilon}(r_{n}-k)k r_{n}}
< \frac{2\tilde{h}_{\ell}L}{\sqrt{\epsilon} L^{2} \tilde{h}_{\ell}} = \frac{2}{\sqrt{\epsilon}L} < 2\sqrt{\epsilon}
\]
So by Lemma \ref{block},
\[
\int \left|\frac{1}{r_{n}+1}\sum_{i=0}^{k}\chi \circ T^{i(r_{n}-k)}\right| d\mu 
\leq \int\left|\frac{1}{L}\sum_{i=0}^{L-1}\chi \circ T^{i(r_{n}-k)p}\right| d\mu + 2\sqrt{\epsilon}
\]
and since $\tilde{h}_{\ell} \leq 2i(r_{n}-k)p < 4L\tilde{h}_{\ell}$ for all $1 \leq i < L$, by Lemma \ref{BH},
\[
\int \left|\frac{1}{r_{n}+1}\sum_{i=0}^{k}\chi \circ T^{i(r_{n}-k)}\right| d\mu < \sqrt{2\epsilon} + 2\sqrt{\epsilon}
\]
Then by Lemma \ref{mix2},
\[
\sum_{q} \sum_{j=0}^{\tilde{h}_{n}-1} \left|\mu(T^{k\tilde{h}_{n}}J_{n,j}^{(q)} \cap A) - \mu J_{n,j}^{(q)} \mu A \right| 
< \sqrt{2\epsilon} + 2\epsilon + \epsilon + 2\left(\sqrt{2\epsilon} + 2\sqrt{\epsilon}\right) + \epsilon + \epsilon < \frac{1}{4}\epsilon_{0}
\]
and by Lemma \ref{mix2inv},
\[
\sum_{q} \sum_{j=0}^{\tilde{h}_{n}-1} \left|\mu(T^{k\tilde{h}_{n}}I_{n,j} \cap A \cap \tilde{C}_{n}^{(q)}) - \mu I_{n,j} \mu (A \cap \tilde{C}_{n}^{(q)}) \right| 
< \frac{1}{4}\epsilon_{0}
\]
Therefore, by Lemma \ref{ts}, for every $\tilde{h}_{n} \leq t < (r_{n}-1)\tilde{h}_{n}$,
\[
\sum_{q} \sum_{j=0}^{\tilde{h}_{n}-1} \left|\mu(T^{t} J_{n,j}^{(q)} \cap A) - \mu J_{n,j}^{(q)} \mu A \right| < 2\left(\frac{1}{4}\epsilon_{0} + \frac{1}{4}\epsilon_{0}\right) = \epsilon_{0}
\]
and likewise
\[
\sum_{q} \sum_{j=0}^{\tilde{h}_{n}-1} \left|\mu(T^{t} I_{n,j} \cap A \cap \tilde{C}_{n}^{(q)}) - \mu I_{n,j} \mu (A \cap \tilde{C}_{n}^{(q)}) \right| < \epsilon_{0}
\]

Finally, consider $(1 - \frac{1}{2}\epsilon_{0})\tilde{h}_{n} \leq t < \tilde{h}_{n}$.  Then by Lemma \ref{mix0all},
\begin{align*}
\sum_{q} \sum_{j=\tilde{h}_{n}-t}^{h_{n}} \left|\mu(T^{t}J_{n,j}^{(q)} \cap A) - \mu J_{n,j}^{(q)} \mu A \right|
&\leq \sum_{q} \sum_{j=r_{n}}^{t} \left|\mu(T^{\tilde{h}_{n}}J_{n,j}^{(q)} \cap A) - \mu J_{n,j}^{(q)} \mu A \right| + \sum_{q} \sum_{j=0}^{r_{n}-1} \mu J_{n,j}^{(q)}\\
&\leq \int \left|\frac{1}{r_{n}+1}\sum_{i=0}^{r_{n}-1}\chi \circ T^{-i}\right| d\mu + \frac{1}{r_{n}+1} + \frac{r_{n}}{\tilde{h}_{n}}
\end{align*}
and therefore
\[
\sum_{q} \sum_{j=0}^{\tilde{h}_{n}-1} \left|\mu(T^{t}J_{n,j}^{(q)} \cap A) - \mu J_{n,j}^{(q)} \mu A \right|
< \epsilon + \frac{h_{n}-t}{h_{n}} + \epsilon + \epsilon
< \epsilon + \frac{1}{2}\epsilon_{0} + 2\epsilon < \epsilon_{0}
\]
Similar reasoning using the second statement in Lemma \ref{mix0all} gives
\[
\sum_{q} \sum_{j=0}^{\tilde{h}_{n}-1} \left|\mu(T^{t}I_{n,j} \cap A \cap \tilde{C}_{n}^{(q)}) - \mu I_{n,j} \mu (A \cap \tilde{C}_{n}^{(q)}) \right|
 < \epsilon_{0} \qedhere
\]
\end{proof}

\section{Partial mixing}

\begin{lemma}\label{mesh2}
Let $T$ be a rank-one transformation such that $B_{n+1}$ through $B_{n+m}$ are constructed using Chacon-like operations with $c_{n+z} = c_{n}$ and $M_{n+z} = M_{n}$ for all $0 \leq z < m$ where $M_{n}$ and $c_{n}$ are such that $M_{n}$ divides $\tilde{h}_{n}$.

For any $1 \leq \ell < \frac{\tilde{h}_{n+m}}{M_{n}}$, 
 there exist disjoint mesh prototypes $J_{n-1,0}^{(q)}$ of $I_{n-1,0}$ for each $-\frac{\tilde{h}_{n}}{2M_{n}} < q \leq \frac{\tilde{h}_{n}}{2M_{n}}$ and a set $Q_{n-1,0}$ which is a union of levels in $C_{n+m}$ such that $T^{j}Q_{n-1,0} \cap Q_{n-1,0} = \varnothing$ for $1 \leq j < \tilde{h}_{n}$ and
\[
T^{\ell M_{n}}I_{n-1,0} = T^{\tilde{h}_{n+m}} Q_{n-1,0} \sqcup \bigsqcup_{q=-\tilde{h}_{n}/2M_{n}+1}^{\tilde{h}_{n}/2M_{n}}  T^{qM_{n}} J_{n-1,0}^{(q)}
\]
If there exists $0 < \beta < 1$ such that $\mu(T^{\ell M_{n}}I_{n,0} \cap I_{n,0}) \leq \beta\mu(I_{n,0})$ and $\mu(T^{(\ell-1) M_{n}}I_{n,0} \cap I_{n,0}) \leq \beta\mu(I_{n,0})$
then $\mu (J_{n-1,0}^{(0)}) \leq \beta\mu(I_{n-1,0})$ and $\mu (J_{n-1,0}^{(-1)}) \leq \beta\mu(I_{n-1,0})$.
\end{lemma}
\begin{proof}
Since $\tilde{h}_{n+1} = (r_{n}+1)\tilde{h}_{n} + s_{n}M_{n} + c_{n+1} - c_{n}$ where $s_{n}$ is the number of subcolumns with $M_{n}$ non-elevation spacers placed over them, and since $c_{n+1} = c_{n}$, it follows that $\tilde{h}_{n+1}$ is divisible by $M_{n}$ as by hypothesis $\tilde{h}_{n}$ is.  The same reasoning gives that $\tilde{h}_{n+z}$ is divisible by $M_{n}$ for all $0 \leq z < m$.

Note that $\ell M_{n} < \tilde{h}_{n+m}$.
Let $S$ be the union of all non-elevation spacers levels placed above $C_{n}$ through $C_{n+m-1}$ and let $S^{\star}$ be the non-elevation spacer levels placed above $C_{n+m}$.
Let
\[
Q_{n,0} = \{ x \in I_{n,0} : \exists~0 < y \leq \ell M_{n}~\text{s.t.}~T^{y}x \in S^{\star} \}
\]
be the set of points in $I_{n,0}$ which $T^{\ell M_{n}}$ pushes through the top of $C_{n+m}$.  For $-\frac{\tilde{h}_{n}}{2M_{n}} < q \leq \frac{\tilde{h}_{n}}{2M_{n}}$, let
\[
Z_{n,0}^{(q)} = \{ x \in I_{n,0} \setminus Q_{n,0} : T^{-qM_{n} + \ell M_{n}}x \in I_{n,0} \}
\]
be the set of points in $I_{n,0}$ such that $T^{\ell M_{n}}$ does not push them through the top of $C_{n+m}$ and leaves them exactly $qM_{n}$ away from $I_{n,0}$.  Since $M_{n}$ divides $\tilde{h}_{n+z}$ for all $0 \leq z < m$,
\[
I_{n,0} \setminus Q_{n,0} \subseteq \bigsqcup_{q=1 - \tilde{h}_{n}/2M_{n}}^{\tilde{h}_{n}/2M_{n}} Z_{n,0}^{(q)}
\]

For $x \in Q_{n,0}$, if $x \in I_{n+m,j}$ then since $T^{y}x \in S^{\star}$ for some $0 < y \leq \ell M_{n} < \tilde{h}_{n+m}$, it follows that the same holds for all $y \in I_{n+m,j}$ since if any point in a level is pushed through the top of the column then all points in the level are (and no point is pushed through twice as $\ell M_{n} < \tilde{h}_{n+m}$).  Therefore $Q_{n,0}$ can be written as union of levels in $C_{n+m}$.  For any $I_{n+m,j} \subseteq Q_{n,0}$, since $T^{\ell M_{n}}$ pushes $I_{n+m,j}$ through the top of $C_{n+m}$ exactly once, it must be that $j \geq \tilde{h}_{n+m} - \ell M_{n}$.  

Set $Q_{n,0}^{\prime} = T^{\ell M_{n} - \tilde{h}_{n+m}}Q$. Then 
 $Q_{n,0}^{\prime}$ is a union of levels in $C_{n+m}$ and $T^{\ell M_{n}}Q_{n,0} = T^{\tilde{h}_{n+m}}Q_{n,0}^{\prime}$.
If $x \in T^{j}Q_{n,0} \cap Q_{n,0}$ for $0 \leq j < \tilde{h}_{n}$ then $x \in I_{n,j} \cap I_{n,0}$ so $j=0$.

By construction, $T^{\ell M_{n}}Z_{n,0}^{(q)} \subseteq T^{qM_{n}}I_{n,0}$.  Set
\[
J_{n,0}^{(q)} = T^{-qM_{n}}\left(T^{\ell M_{n}}Z_{n,0}^{(q)}\right) \subseteq I_{n,0}
\]
and set, writing $s_{n-1,y}$ for the non-elevation spacer amounts added above $C_{n-1}$,
\[
J_{n-1,0}^{(q)} = \bigsqcup_{i=0}^{r_{n}-1} T^{i\tilde{h}_{n-1}+\sum_{y=0}^{i-1}s_{n-1,y}}J_{n,0}^{(q)}
\]
Then $J_{n-1,0}^{(q)}$ is a mesh prototype for $I_{n,0}$ (since $I_{n-1,0}^{[i]} = T^{i\tilde{h}_{n-1}+\sum_{y=0}^{i-1}s_{n-1,y}}I_{n,0}$) and
\begin{align*}
T^{\ell M_{n}}I_{n-1,0} &= T^{\ell M_{n}}\left(\bigsqcup_{i=0}^{r_{n-1}} T^{i\tilde{h}_{n-1}+\sum_{y=0}^{i-1}s_{n-1,y}}I_{n,0}\right) \\
&= \bigsqcup_{i=0}^{r_{n-1}} T^{i\tilde{h}_{n-1}+\sum_{y=0}^{i-1}s_{n-1,y}} \left( T^{\tilde{h}_{n+m}}Q_{n,0}^{\prime} \sqcup \bigsqcup_{q} T^{qM_{n}} J_{n,0}^{(q)}\right) \\
&= T^{h_{n+m}}\left(\bigsqcup_{i=0}^{r_{n-1}}  T^{i\tilde{h}_{n-1}+\sum_{y=0}^{i-1}s_{n-1,y}} Q_{n,0}^{\prime}\right)
\sqcup \bigsqcup_{q} T^{qM_{n}}J_{n-1,0}^{(q)}
\end{align*}
Set $Q_{n-1,0} = \bigsqcup_{i=0}^{r_{n-1}}  T^{i\tilde{h}_{n-1}+\sum_{y=0}^{i-1}s_{n-1,y}} Q_{n,0}^{\prime}$.  Then $Q_{n-1,0}$ is a union of levels in $C_{n+m}$ since $i\tilde{h}_{n-1}+\sum_{y=0}^{i-1}s_{n-1,y} < \tilde{h}_{n}$ for all $0 \leq i < r_{n-1}$.

Since for $q < q^{\prime}$, $J_{n,0}^{(q)} \cap J_{n,0}^{(q^{\prime})} \subseteq T^{\ell M_{n}+qM_{n}}(Z_{n,0}^{(q)} \cap T^{(q^{\prime} - q)M_{n}}Z_{n,0}^{(q^{\prime})})$, if the intersection were nonempty then $(q^{\prime} - q)M_{n}$ would have to be a multiple of $\tilde{h}_{n}$ as $Z_{n,0}^{(q)} \subseteq I_{n,0}$.  But that is impossible so it follows that the $J_{n,0}^{(q)}$ and hence the $J_{n-1,0}^{(q)}$ are mutually disjoint.

Assume now that
there exists $0 < \beta < 1$ such that $\mu(T^{\ell M_{n}}I_{n,0} \cap I_{n,0}) \leq \beta\mu(I_{n,0})$ and $\mu(T^{(\ell-1) M_{n}}I_{n,0} \cap I_{n,0}) \leq \beta\mu(I_{n,0})$.  Since $J_{n,0}^{(0)} \subseteq I_{n,0} \cap T^{\ell M_{n}}I_{n,0}$, it follows that $\mu(J_{n,0}^{(0)}) = \mu(J_{n,0}^{(0)} \cap I_{n,0} \cap T^{\ell M_{n}}I_{n,0}) \leq \beta \mu(I_{n,0})$.  Likewise, since $J_{n,0}^{(-1)} \subseteq I_{n,0} \cap T^{(\ell - 1)M_{n}}I_{n,0}$, the same holds for $J_{n,0}^{(-1)}$.
\end{proof}

\begin{proposition}\label{pmix1}
Let $T$ be an elevated staircase-Chacon-like hybrid transformation and $A$ and $D$ be unions of levels in some column $C_{n}$ and let $\epsilon > 0$.  If for some $m > 2$,
\begin{itemize}
\item $A$ and $D$ are uniformly mixed within $\epsilon$ at stages $n-2$ and $n-1$
\item $A$ is uniformly mixed within $\epsilon$ at stage $n+m$
\item $B_{n+1}$ through $B_{n+m}$ are constructed using Chacon-like operations
\item $M_{n+z} = 2\tilde{h}_{n-2}$ for all $0 \leq z < m$
\item for all $0 \leq z < m$, $c_{n+z} = c_{n}$ is the unique minimal positive integer such that $c_{n} \geq c_{n-1} + r_{n-1}$ and $2\tilde{h}_{n-2}$ divides $h_{n} + c_{n} = \tilde{h}_{n}$
\item there exists $0 < \beta < 1$ such that $\mu(T^{\ell M_{n}}I_{n,0} \cap I_{n,0}) \leq \beta\mu(I_{n,0})$ for all $\frac{\tilde{h}_{n}}{M_{n}} - 1 \leq \ell < \frac{\tilde{h}_{n+m}}{M_{n}}$
\end{itemize}
then for all $(1 - \frac{1}{2}\epsilon)\tilde{h}_{n-2} \leq t < (1 - \frac{1}{2}\epsilon)\tilde{h}_{n+m}$,
\[
\mu(T^{t}D \cap A) \geq (1 - \beta)\mu D \mu A - 5\epsilon
\]
\end{proposition}
\begin{proof}
The claim is immediate for $(1 - \frac{1}{2}\epsilon)\tilde{h}_{n-2} \leq t < (1 - \frac{1}{2}\epsilon)\tilde{h}_{n}$ since
\[
\left|\mu(T^{t}D \cap A) - \mu D \mu A\right|
\leq \sum_{j=0}^{h_{n}-1}\left|\mu(T^{t}I_{n,j} \cap A) - \mu I_{n,j} \mu A\right| < \epsilon
\]
so $\mu(T^{t}D \cap A) \geq \mu D \mu A - \epsilon$.

Let $\tilde{h}_{n} \leq t < (1 - \frac{1}{2}\epsilon)\tilde{h}_{n+m}$.  Let $\ell$ be the unique integer such that $\ell M_{n} \leq t < (\ell + 1)M_{n}$.  Then $\frac{\tilde{h}_{n}}{M_{n}} \leq \ell < \frac{\tilde{h}_{n+m}}{M_{n}}$ and let $y = t - \ell M_{n}$ be such that $0 \leq y < M_{n}$.  Let $J_{n-1,0}^{(q)}$ be the mesh prototypes and $Q_{n-1,0}$ be the set from Lemma \ref{mesh2}.  Then
\[
T^{t}I_{n-1,0} = T^{y}T^{\ell M_{n}}I_{n-1,0} = T^{\tilde{h}_{n+m} + y}Q_{n-1,0} \sqcup \bigsqcup_{q=1-\tilde{h}_{n}/2M_{n}}^{\tilde{h}_{n}/2M_{n}} T^{qM_{n}+y} J_{n-1,0}^{(q)}
\]
Let $\tilde{C}_{n-1}^{(q)} = \bigsqcup_{j=0}^{\tilde{h}_{n-1}-1} J_{n-1,j}^{(q)}$ be the disjoint meshes from the prototypes.
By Lemma \ref{meshdrop}, $C_{n-2} \cap \tilde{C}_{n-1}^{(q)}$ is a mesh of $C_{n-2}$ which we will write as $\tilde{C}_{n-2}^{(q)} = \bigsqcup_{j=0}^{\tilde{h}_{n-2}-1}J_{n-2,j}^{(q)}$.  

Since $Q_{n-1,0}$ is a union of levels in $C_{n+m}$ and $Q_{n-1,0} \cap T^{j}Q_{n-1,0} = \varnothing$ for $1 \leq j < \tilde{h}_{n}$, set $Q_{n-1,j} = T^{j}Q_{n-1,0}$ for each $0 \leq j < \tilde{h}_{n-1}$ and $Q_{n-1} = \bigsqcup_{j=0}^{\tilde{h}_{n-1}} Q_{n-1,j}$ is a union of levels in $C_{n+m}$.

Since $D$ can be written as a union of levels in $C_{n-1}$, there exists $\mathcal{D} \subseteq \{ 0, 1, \ldots, \tilde{h}_{n-1}-1 \}$ such that $D = \bigsqcup_{j \in \mathcal{D}} I_{n-1,j}$.  Then
\begin{align*}
T^{t}D = \bigsqcup_{j\in\mathcal{D}} T^{t}I_{n,j}
&= \bigsqcup_{j\in\mathcal{D}}\left(T^{\tilde{h}_{n+m}+y}Q_{n-1,j} \sqcup \bigsqcup_{q} T^{qM_{n}+y}J_{n-1,j}^{(q)}\right) \\
&= T^{\tilde{h}_{n+m}+y}\left(D \cap Q_{n-1}\right) \sqcup \bigsqcup_{q} T^{qM_{n}+y}\left(D \cap \tilde{C}_{n-1}^{(q)}\right)
\end{align*}
and therefore, since $D \subseteq C_{n-2}$,
\begin{align*}
\mu(T^{t}D \cap A) \geq
\mu(T^{\tilde{h}_{n+m}+y}(D \cap Q_{n-1}) \cap A) 
&+ \sum_{q : \tilde{h}_{n-2} \leq |qM_{n}+y| < \tilde{h}_{n-1}} \mu(T^{qM_{n}+y}(D \cap \tilde{C}_{n-2}^{(q)}) \cap A) \\
&+ \sum_{q : |qM_{n} + y| \geq \tilde{h}_{n-1}} \mu(T^{qM_{n}+y}(D \cap \tilde{C}_{n-1}^{(q)}) \cap A)
\end{align*}

Since $Q_{n-1,j}$ are disjoint sets which are unions of levels in $C_{n+m}$, then
\begin{align*}
&\left|\mu(T^{\tilde{h}_{n+m}+y}(D \cap Q_{n-1}) \cap A) - \mu(D \cap Q_{n-1})\mu A\right|
\leq
\sum_{j=0}^{\tilde{h}_{n-1}-1}\left|\mu(T^{\tilde{h}_{n+m}+y}Q_{n-1,j} \cap A) - \mu Q_{n-1,j}\mu A\right|  \\
&\quad\quad\quad\quad= 
\sum_{j=0}^{\tilde{h}_{n-1}-1}\left|\sum_{j^{\prime} : I_{n+m,j^{\prime}}\subseteq Q_{n-1,j}} \mu(T^{\tilde{h}_{n+m}+y}I_{n+m,j^{\prime}} \cap A) - \mu I_{n+m,j^{\prime}}\mu A\right|  \\
&\quad\quad\quad\quad\leq \sum_{j=0}^{\tilde{h}_{n+m}-1}\left|\mu(T^{\tilde{h}_{n+m}+y}I_{n+m,j} \cap A) - \mu I_{n+m,j} \mu A \right| < \epsilon
\end{align*}
since $A$ is uniformly mixed at stage $n+m$ and $\tilde{h}_{n+m} \leq \tilde{h}_{n+m} + y < \tilde{h}_{n+m} + M_{n} < (1 - \epsilon)\tilde{h}_{n+m+1}$.

By Lemma \ref{ineq}, since $A$ is uniformly mixed within $\epsilon$ at stage $n-2$ (and at stage $n-1$ for when $qM_{n}+y \geq (1 - \frac{\epsilon}{2})\tilde{h}_{n-2}$),
\[
\sum_{q : \tilde{h}_{n-2} \leq qM_{n}+y < \tilde{h}_{n-1}} \left|\mu(T^{qM_{n}+y}(D \cap \tilde{C}_{n-2}^{(q)}) \cap A) - \mu(D \cap \tilde{C}_{n-2}^{(q)})\mu A\right| < \epsilon
\]
and since $A$ is uniformly mixed at stage $n-1$ and $|qM_{n} + y| \leq \frac{1}{2}\tilde{h_{n}} + M_{n} = \frac{1}{2}\tilde{h}_{n} + 2\tilde{h}_{n-2} < (1 - \frac{1}{2}\epsilon)\tilde{h}_{n}$,
\[
\sum_{q : \tilde{h}_{n-1} \leq qM_{n}+y} \left|\mu(T^{qM_{n}+y}(D \cap \tilde{C}_{n-1}^{(q)}) \cap A) - \mu(D \cap \tilde{C}_{n-1}^{(q)})\mu A\right| < \epsilon
\]
Since $D$ is uniformly mixed at stages $n-2$ and $n-1$, using the second condition of uniformly mixed within $\epsilon$,
\begin{align*}
&\sum_{q : \tilde{h}_{n-2} \leq -(qM_{n}+y) < \tilde{h}_{n-1}} \left|\mu(T^{qM_{n}+y}(D \cap \tilde{C}_{n-2}^{(q)}) \cap A) - \mu(D \cap \tilde{C}_{n-2}^{(q)})\mu A\right| \\
&\quad\quad = \sum_{q : \tilde{h}_{n-2} \leq -(qM_{n}+y) < \tilde{h}_{n-1}} \left|\mu(D \cap \tilde{C}_{n-2}^{(q)} \cap T^{-(qM_{n}+y)}A) - \mu(D \cap \tilde{C}_{n-2}^{(q)})\mu A\right|
< \epsilon
\end{align*}
and likewise since $D$ is uniformly mixed at stage $n-1$,
\[
\sum_{q : \tilde{h}_{n-1} \leq -(qM_{n}+y)} \left|\mu(T^{qM_{n}+y}(D \cap \tilde{C}_{n-1}^{(q)}) \cap A) - \mu(D \cap \tilde{C}_{n-1}^{(q)})\mu A\right| < \epsilon
\]
Since $y < M_{n}$, for every $q \ne 0,-1$ it holds that $|qM_{n}+y| \geq M_{n}$ and also that one of $|-M_{n}+y| \geq \frac{1}{2}M_{n}$ or $|y| \geq \frac{1}{2}M_{n}$.  Therefore $\sum_{q : |qM_{n} + y| < \tilde{h}_{n-2}} \mu(D \cap \tilde{C}_{n-1}^{(q)}) \leq \mu(D \cap \tilde{C}_{n-1}^{(q^{\prime})})$ where $q^{\prime} \in \{ 0, -1 \}$.  From Lemma \ref{mesh2}, $\mu J_{n-1,0}^{(q^{\prime})} \leq \beta\mu I_{n-1,0}$ and so $\mu J_{n-1,j}^{(q^{\prime})} \leq \beta\mu I_{n-1,j}$ for all $0 \leq j < \tilde{h}_{n-1}-1$.  Therefore, as $D$ is a union of levels,
\[
\mu\left(D \cap \tilde{C}_{n-1}^{(q^{\prime})}\right) \leq \beta\mu D
\]

Since $T^{qM_{n}}I_{n,0} = T^{\tilde{h}_{n+m}}Q_{n-1,0} \sqcup \bigsqcup_{q} T^{qM_{n}}J_{n-1,0}^{(q)}$ and $T$ is measure-preserving, $\mu I_{n,0} = \mu Q_{n-1,0} + \sum_{q} \mu J_{n-1,0}^{(q)}$.
 and therefore
\[
\mu(D) = \mu\left(D \cap \left(Q_{n-1} \sqcup \bigsqcup_{q} \tilde{C}_{n-1}^{(q)}\right)\right)
= \mu\left(D \cap \left(Q_{n-1} \sqcup \bigsqcup_{q\ne q^{\prime}} \tilde{C}_{n-1}^{(q)}\right)\right) + \mu\left(D \cap \tilde{C}_{n-1}^{(q^{\prime})}\right)
\]
meaning that
\[
\mu\left(D \cap \left(Q_{n-1} \sqcup \bigsqcup_{q\ne q^{\prime}} \tilde{C}_{n-1}^{(q)}\right)\right) 
= \mu(D) - \mu\left(D \cap \tilde{C}_{n-1}^{(q^{\prime})}\right)
\geq \mu D - \beta\mu D
\]

Putting all of this together (and using that $D\cap\tilde{C}_{n-2}^{(q)} = D\cap\tilde{C}_{n-1}^{(q)}$), since the meshes are disjoint and disjoint from $Q_{n-1}$,
\begin{align*}
\mu(T^{t}D \cap A) &\geq
\mu(D \cap Q_{n-1})  \mu A 
+ \sum_{q \ne q^{\prime}} \mu(D \cap \tilde{C}_{n-1}^{(q)}) \mu A - 5\epsilon \\
&= \mu\left(D \cap \left(Q_{n-1} \sqcup \bigsqcup_{q\ne q^{\prime}} \tilde{C}_{n-1}^{(q)}\right)\right) \mu A - 5\epsilon 
\geq (1 - \beta) \mu D \mu A - 5\epsilon \qedhere
\end{align*}
\end{proof}

\section{Low complexity partial mixing}

Fix an integer $\gamma \geq 2$.  We will construct an elevated staircase-Chacon-like hybrid transformation $T_{\gamma}$.

Let $B_{0} = 0$ be the starting word.  Let $S_{0} = -1$ and $Q_{0} = 0$.  Set $r_{-1} = 3$ and $c_{0} = 1$ and $h_{-1} = 3$.

Assume that the transformation has been constructed up to $B_{Q_{n}}$ and that $c_{Q_{n}}$ has been specified.  Construct $B_{Q_{n}+1}$ and $B_{Q_{n}+2}$ using the staircase operation with $r_{Q_{n}} = r_{Q_{n}+1} = \tilde{h}_{S_{n}}$ and $c_{Q_{n}+1} = c_{Q_{n}} + r_{Q_{n}}$.  Set $S_{n+1} = Q_{n} + 2$.

Let $M_{S_{n+1}} = 2\tilde{h}_{Q_{n}}$.  Let $c_{S_{n+1}}$ be the minimal integer such that $c_{S_{n+1}} > c_{S_{n+1}-1} + r_{S_{n+1}-1}$ and $M_{S_{n+1}}$ divides $h_{S_{n+1}} + c_{S_{n+1}}$.
By Proposition \ref{3prime}, $p(\tilde{h}_{S_{n+1}})$
 is now completely determined.  Choose $y_{n+1} > S_{n+1} + 2$ such that 
\[
3^{S_{n+1} - y_{n+1}}  \frac{p(\tilde{h}_{S_{n+1}})}{\tilde{h}_{S_{n+1}}} < \frac{1}{n+1}
\]

Construct $B_{t}$ for $S_{n} < t \leq y_{n+1}$ using the $\gamma$-Chacon-like operation with $M_{t} = M_{S_{n+1}}$ and $c_{t} = c_{S_{n+1}}$.
Choose $Q_{n+1} > y_{n+1} + n+1$ such that $\tilde{h}_{S_{n+1}}^{2} < \frac{1}{n+1} 3^{Q_{n+1}}$.
Construct $B_{t}$ for $y_{n+t} < t \leq Q_{n+1}$ using the $\gamma$-Chacon-like operation with $M_{t} = M_{S_{n+1}}$ and $c_{t} = c_{S_{n+1}}$.
Set $c_{Q_{n+1}} = |B_{y_{n+1}}| + 1$.

Proceed inductively for all $n$.

\subsection{Basic properties}

\begin{lemma}\label{good}
$T_{\gamma}$ satisfies
\[
\frac{|B_{y_{n+1}}|+1}{|B_{Q_{n+1}}|} < 3^{-n-1}
\quad\quad
\text{and}
\quad\quad
\frac{r_{Q_{n+1}}^{2}}{|B_{Q_{n+1}}|} < \frac{1}{n+1}
\quad\quad
\text{and}
\quad\quad
\frac{r_{Q_{n+1}+1}^{2}}{|B_{Q_{n+1}+1}|} < \frac{1}{n+1}
\]
\end{lemma}
\begin{proof}
By construction, since $|B_{t+1}| \geq 3|B_{t}|$ for all $t$ and $Q_{n+1} \geq y_{n+1} + n + 1$,
\[
\frac{|B_{y_{n+1}}|+1}{|B_{Q_{n+1}}|} < 3^{-Q_{n+1}+y_{n+1}+1} < 3^{-n}
\quad\quad
\text{and}
\quad\quad
\frac{\tilde{h}_{S_{n}}^{2}}{|B_{Q_{n+1}}|} < \frac{\frac{1}{n+1}3^{Q_{n+1}}}{3^{Q_{n+1}}} = \frac{1}{n+1}
\]
and we have chosen $r_{Q_{n+1}} = r_{Q_{n+1}+1} = \tilde{h}_{S_{n}}$.
\end{proof}

\begin{lemma}\label{Moverh}
$T_{\gamma}$ satisfies $\frac{M_{S_{n}}}{h_{S_{n}-1}} \leq 3^{-n}$.
\end{lemma}
\begin{proof}
Since $|B_{t+1}| > 3|B_{t}|$ (and $c_{n+1} > c_{n}$),
\[
h_{S_{n}-1} > r_{S_{n}-2}\tilde{h}_{S_{n}-2} = r_{Q_{n-1}}~\frac{1}{2}M_{S_{n}} = h_{S_{n-1}}~\frac{1}{2}M_{S_{n}}
\geq 3^{S_{n-1}}~\frac{1}{2}M_{S_{n}}
\]
and $S_{n-1} \geq n+1$ so the claim follows.
\end{proof}

\begin{lemma}\label{cr}
$T_{\gamma}$ satisfies
$
\sum \frac{c_{n}}{|B_{n}|} < \infty
$.
\end{lemma}
\begin{proof}
By Lemma \ref{good},
$
\frac{c_{Q_{n}}}{|B_{Q_{n}}|} = \frac{|B_{y_{n}}|+1}{|B_{Q_{n}}|} < 3^{-n}
$.
By construction, using that $|B_{t}| > 3^{t}$ and $Q_{n} \geq n$,
\[
\frac{c_{Q_{n}+1}}{|B_{Q_{n}+1}|} 
= \frac{c_{Q_{n}} + r_{Q_{n}}}{|B_{Q_{n}+1}|}
< \frac{c_{Q_{n}}+r_{Q_{n}}}{r_{Q_{n}}|B_{Q_{n}}|}
= \frac{1}{r_{Q_{n}}}\frac{c_{Q_{n}}}{|B_{Q_{n}}|} + \frac{1}{|B_{Q_{n}}|}
< 3^{-n-1} + 3^{-n-1} < 3^{-n}
\]
and similarly, by Lemma \ref{Moverh},
\[
\frac{c_{S_{n+1}}}{|B_{S_{n+1}}|}
< \frac{c_{Q_{n}+1} + r_{Q_{n}+1} + M_{S_{n+1}}}{|B_{S_{n+1}}|}
< \frac{c_{Q_{n}+1} + r_{Q_{n}+1} + M_{S_{n+1}}}{r_{Q_{ n}+1}|B_{Q_{n}+1}|}
< 3^{-n-1} + 3^{-Q_{n}-1} + 3^{-n-1} < 3^{-n}
\]
For $0 < t < Q_{n+1} - S_{n+1}$,
$
\frac{c_{S_{n+1}+t}}{|B_{S_{n+1}+t}|}
< \frac{c_{S_{n+1}}}{3^{t}|B_{S_{n+1}}|}
< 3^{-n-t}
$.
Therefore
\[
\sum_{t=0}^{\infty} \frac{c_{t}}{|B_{t}|}
= \sum_{n=0}^{\infty} \sum_{t=Q_{n}}^{Q_{n+1}-1} \frac{c_{t}}{|B_{t}|}
< \sum_{n=0}^{\infty} \left(3 \cdot 3^{-n} + \sum_{t=1}^{Q_{n+1}-S_{n+1}-1} 3^{-n-t}\right)
< \sum_{n=0}^{\infty} \left(3^{-n+1} + 3^{-n}\frac{1}{2}\right) < \infty \qedhere
\]
\end{proof}

\begin{proposition}
$T_{\gamma}$ is on a finite measure space.
\end{proposition}
\begin{proof}
The total number of spacers added at each stage is $c_{t}$ levels of the same size as those in $B_{t}$ and either approximately $\frac{1}{2}r_{t}^{2}$ levels the size of those in $B_{t+1}$ or $M_{t}$ levels the size of those in $B_{t+1}$.

The measure of the elevation levels is then $\sum_{t} \frac{c_{t}}{h_{t}}$ which is finite by Lemma \ref{cr}.  By Lemma \ref{good}, $\frac{r_{n}^{2}}{h_{n}} \to 0$ for those $n$ for which $B_{n+1}$ is constructed using the staircase operation so the spacer levels in the staircase operations have finite total measure.

At each Chacon-like form operation, the non-elevation levels have measure at most $\frac{M_{t}}{|B_{t}|}$ since at most $M_{t}$ spacers are added above each subcolumn.  By Lemma \ref{Moverh}, $\sum \frac{M_{S_{n}}}{|B_{S_{n}}|} < \infty$.  Since $|B_{t+1}| > 3|B_{t}|$, the total measure added by the non-elevation spacers in Chacon-like form operations is, as $M_{t} = M_{S_{n+1}}$ for $S_{n+1} \leq t < Q_{n+1}$,
\[
\sum_{n=0}^{\infty} \sum_{t=S_{n+1}}^{Q_{n+1}-1} \frac{M_{t}}{|B_{t}|}
\leq \sum_{n=0}^{\infty} \sum_{t=0}^{Q_{n+1}-S_{n+1}-1} \frac{M_{S_{n+1}}}{|B_{S_{n+1}}|} \frac{|B_{S_{n+1}}|}{|B_{S_{n+1}+t}|}
\leq \sum_{n=0}^{\infty} \frac{M_{S_{n+1}}}{|B_{S_{n+1}}|} \sum_{t=0}^{\infty} 3^{-t} < \infty \qedhere
\]
\end{proof}

\subsection{Word complexity}

\begin{proposition}\label{comp}
$T_{\gamma}$ satisfies
$
\liminf \nicefrac{p(q)}{q} \leq 2 + \nicefrac{1}{(\gamma + 1)}
$
\end{proposition}
\begin{proof}
By Proposition \ref{2full}, since $|B_{y_{n}}| < c_{Q_{n}}$,
\[
\frac{p\left(\frac{1}{2}\tilde{h}_{n+1}\right)}{\frac{1}{2}\tilde{h}_{n+1}} < 2 + \frac{1}{\gamma + 1} + 3^{S_{n}-y_{n}}\frac{p(|B_{S_{n}}|)}{|B_{S_{n}}|}
\]
so by our choice of $y_{n}$, then
\[
\frac{p\left(\frac{1}{2}\tilde{h}_{n+1}\right)}{\frac{1}{2}\tilde{h}_{n+1}}   \leq 2 + \frac{1}{\gamma + 1} + \frac{1}{n} \to 2 + \frac{1}{\gamma + 1} \qedhere
\]
\end{proof}


\begin{remark}
The $T_{\gamma}$ also satisfy $\limsup \nicefrac{p(q)}{q} = \infty$ which is the primary motivation for Conjecture \ref{conj2}.

The staircase construction, which necessarily occurs infinitely often in our constructions with cut parameter tending to infinity in order to ensure partial mixing, forces there to be a sequence of times along which $p(q)/q \to \infty$.  Specifically, at each staircase stage, $B_{n}1^{c_{n}+i}$ is a right-special word for each $0 \leq i < r_{n}$ and therefore for $c_{n} + r_{n} \leq q < |B_{n}|$ we effectively have $p(q+1) - p(q) \geq r_{n} -1$ (we omit the details). As $|B_{n}| \to \infty$ rapidly, this forces $p(|B_{n}|) / |B_{n}| \to \infty$ for $n$ such that a staircase construction is about to be implemented.

In essence, the tension here is that the Chacon-like stages are required to reduce the word complexity to as close to $2$ as possible while the staircase constructions are required to ensure partial mixing. Conjecture \ref{conj1} in essence states that this sort of tension should be unavoidable regardless of what subshift is exhibited since the partial mixing condition should effectively mean that ``mixing occurs infinitely often" and, as shown in \cite{creutzmixing}, outright mixing implies $\lim p(q)/q = \infty$.
\end{remark}

\subsection{Partial mixing}

\begin{notation}
For two words $u$ and $v$ of the same length, $d_{H}(u,v) = |\{ i : u_{[i,i+1)} \ne v_{[i,i+1)} \}|$ is the Hamming distance.
\end{notation}

\begin{proposition}\label{shifts}
Let $P \geq 1$ and let $T$ be 
the subshift on the alphabet $\{0, 1, \ldots, P \}$ constructed using $B_{0} = 023\ldots P$ and $B_{n+1} = (B_{n}1)^{\gamma}B_{n}(B_{n}1)^{\gamma}B_{n}$.
Let $W_{n} = (B_{n}1)^{\gamma}B_{n}$ and for $1 \leq i < |W_{n}|$, let $W_{n,i} = (W_{n})_{[i,|W|_{n}]}(W_{n})_{[1,i)}$ be the word obtained by shifting $W_{n}$ to the left $i$ places (with wraparound).

Let $W$ be the right-infinite word which has $W_{n}$ at the origin for every $n$.  
Then there exists a constant $\kappa_{\gamma} > 0$ depending only on $\gamma$ (and not on $P$) such that if $i \geq 1$ and  $n$ is such that $|W_{n}| > i$ then
$d_{H}(W_{[1,|W_{n}|]},(\sigma^{i}W)_{[1,|W_{n}|]}) \geq \kappa_{\gamma}|W_{n}|$.
\end{proposition}
\begin{proof}
Let $|u|_{0}$ be the number of $0$s in a word $u$ and likewise $|u|_{1}$ be the number of $1$s.
Since $|B_{0}|_{0} = 1$ and $|B_{n+1}|_{0} = (2\gamma+2)|B_{n}|_{0}$, then $|B_{n}|_{0} = (2\gamma+2)^{n}$.  Since $|B_{0}|_{1} = 0$ and $|B_{n+1}|_{1} = 2\gamma + (2\gamma+2)|B_{n}|_{1}$, then $|B_{n}|_{1} = 2\gamma(1 + (2\gamma+2) + \cdots + (2\gamma+2)^{n-1}) = 2\gamma\frac{(2\gamma+2)^{n}-1}{2\gamma+1}$.  Then for $n \geq 1$,
\[
\frac{|B_{n}|_{1}}{|B_{n}|_{0}} = \frac{2\gamma - \frac{2\gamma}{(2\gamma+2)^{n-1}}}{2\gamma + 1} 
> \frac{2\gamma - 1}{2\gamma+1}
\quad\quad\text{and}\quad\quad
\frac{P|B_{n}|_{0}}{|B_{n}|} =
\frac{1}{1 + \frac{|B_{n}|_{1}}{P|B_{n}|_{0}}}
> \frac{1}{1 + \frac{2\gamma}{P(2\gamma+1)}}
\]

For a word $u$ which is a concatenation of 1s and $B_{n}$ blocks, since $0$ is always followed by $23\ldots P$ in such a block, the number of letters in $u$ which are not $1$ is $P|u|_{0}$.
For words $u$ and $v$ of the same length, let $d_{0}(u,v) = \frac{1}{P|u|_{0}}|\{ i : u_{[i,i+1)} \ne 1, v_{[i,i+1)} \ne u_{[i,i+1)} \}|$.

Since $B_{n}$ begins with $0$ and $11$ never appears, every occurrence of a $1$ in $1B_{n}$ coincides with an occurrence of $0$ in $B_{n}1$.  As there are then at least $P$ discrepancies between the words starting at that location,
\[
d_{0}(B_{n}1,1B_{n}) \geq \frac{P|B_{n}|_{1}}{P|B_{n}|_{0}} > \frac{2\gamma - 1}{2\gamma+1}
\]
and likewise $d_{0}(1B_{n},B_{n}1) > \frac{2\gamma - 1}{2\gamma+1}$.

We claim now that for every $n$ and every word $w$ in the language such that $|w| = |B_{n}|$, if $w \ne B_{n}$ then $d_{0}(B_{n},w) \geq \frac{1}{2\gamma+2}~\frac{2\gamma-1}{2\gamma+1}$.  When $n=0$, the claim is immediate since $d_{0}(a,1) = 1$ for any letter $a \ne 1$.

Assume the claim holds for $n$.  Let $w$ be a word in the language such that $|w| = |B_{n+1}|$ and $w \ne B_{n+1}$.  Then there exist words $w_{j}$ for $1 \leq j \leq 2\gamma + 2$ and $z_{j}$ for $1 \leq j \leq \gamma$ and $\gamma + 2 \leq j \leq 2\gamma+1$ with $|w_{j}| = |B_{n}|$ and $|z_{j}| = 1$ such that
\[
w = \left(\prod_{j=1}^{\gamma}w_{j}z_{j}\right)w_{\gamma+1} \left(\prod_{j=\gamma+2}^{2\gamma+1}w_{j}z_{j}\right)w_{2\gamma+2}
\]
Since $B_{n+1} = (B_{n}1)^{\gamma}B_{n}(B_{n}1)^{\gamma}B_{n}$ and $d_{0}$ ignores the locations where $B_{n+1}$ reads $1$,
\[
d_{0}(B_{n+1},w) = \frac{|\{ i : (B_{n+1})_{[i,i+1)} \ne 1, w_{[i,i+1)} \ne , i~\text{is in a $B_{n}$ in $B_{n+1}$} \}|}{(2\gamma+2)P|B_{n}|_{0}}
= \frac{1}{2\gamma+2}\sum_{j=1}^{2\gamma+2} d_{0}(B_{n},w_{j})
\]
Suppose $w_{j} \ne B_{n}$ for all $j$.  Then by the inductive assumption, $d_{0}(B_{n+1},w) \geq \frac{1}{2\gamma+2}\sum_{j=1}^{2\gamma+2} \frac{1}{2\gamma+2}~\frac{2\gamma-1}{2\gamma+1} = \frac{1}{2\gamma+2}~\frac{2\gamma-1}{2\gamma+1}$.  So we may assume that $w_{j} = B_{n}$ for some $j$.
Take $j$ to be the minimal location at which $w_{j} = B_{n}$.

Consider first the case when $j \ne 1, j \ne \gamma+1, j \ne 2\gamma+2$.  Since $B_{n}$ is always preceded by either $B_{n}$ or $B_{n}1$, and since $w_{j-1} \ne B_{n}$, it must be that $w_{j-1}z_{j}w_{j} = aB_{n}B_{n}$ for some letter $a$.  Since $B_{n}B_{n}B_{n}$ never appears, then $w_{j-1}z_{j}w_{j} = 1B_{n}B_{n}$ so $w$ and $B_{n+1}$ have at least one location where $w$ reads $1B_{n}$ and $B_{n+1}$ reads $B_{n}1$.

Now consider when $j = 1$.  Let $k$ minimal such that $w_{k} \ne B_{n}$.  Then $w$ has prefix $(B_{n}1)^{k-2}B_{n}$ but does not have prefix $(B_{n}1)^{k-1}$ so, as $B_{n}$ is always followed by $B_{n}$ or $1B_{n}$, then $w$ has prefix $(B_{n}1)^{k-2}B_{n}B_{n}$ hence $B_{n}B_{n}1$ as $B_{n}B_{n}$ is always followed by $1$.
If $k = \gamma+2$ then $w$ would have prefix $(B_{n}1)^{\gamma}B_{n}B_{n}1$ but then $w_{k} = B_{n}$ so $k \ne \gamma + 2$.  Therefore $z_{k-1}$ exists and $z_{k-1}w_{k} = B_{n}1$.  So there is at least one location where $w$ reads $B_{n}1$ and $B_{n+1}$ reads $1B_{n}$.

Then in both cases,
\[
d_{0}(B_{n+1},w) \geq \frac{1}{2\gamma+2}d_{0}(B_{n}1,1B_{n}) \geq \frac{1}{2\gamma+2}~\frac{2\gamma-1}{2\gamma+1}
\]
By induction this holds for all $n$ and all words $w \ne B_{n}$ with $|w| = |B_{n}|$.

Since $B_{n+1} = W_{n}W_{n}$, for all words $w$ with $|w| = |W_{n}|$, then $d_{0}(W_{n},w) \geq \frac{1}{2\gamma+2}~\frac{2\gamma}{2\gamma+1}$ and $\frac{P|W_{n}|_{0}}{|W_{n}|} = \frac{P|B_{n+1}|_{0}}{|B_{n+1}|}$ and therefore
\begin{align*}
\frac{d_{H}(W_{n},W_{n,i})}{|W_{n}|} &\geq d_{0}(W_{n},W_{n,i})\frac{P|W_{n}|_{0}}{|W_{n}|}
> \frac{1}{2\gamma+2}~\frac{2\gamma-1}{2\gamma+1}~\frac{1}{1 + \frac{2\gamma}{P(2\gamma+1)}} \\
&\geq \frac{1}{2\gamma+2}~\frac{2\gamma-1}{2\gamma+1}~\frac{1}{1 + \frac{2\gamma}{2\gamma+1}} 
= \frac{2\gamma-1}{(2\gamma+2)(4\gamma+1)}
\end{align*}
Let $W$ be the right-infinite word that is the limit of $W_{n}$.  Then for any $i \geq 1$, taking $n$ sufficiently large that $|W_{n}| > i$, $\sigma^{i}W$ has $W_{n,i}$ at the origin.  Then for such $n$,
\[
d_{H}(W_{[1,|W_{n}|]},(\sigma^{i}W)_{[1,|W_{n}|]}
= d_{H}(W_{n},W_{n,i}) >  \frac{2\gamma-1}{(2\gamma+2)(4\gamma+1)}|W_{n}|
\]
so the claim follows setting $\kappa_{\gamma} =  \frac{2\gamma-1}{(2\gamma+2)(4\gamma+1)}$.
\end{proof}


\begin{remark}
The above transformations, in the case $P=1$ so on a two letter alphabet, are precisely the examples exhibited in \cite{creutzrankone} demonstrating rank-one weak mixing transformations with word complexity $\liminf \nicefrac{p(q)}{q} = 1$ and $\limsup \nicefrac{p(q)}{q} < 1.5 + \epsilon$.
\end{remark}

\begin{lemma}\label{norigid}
$T_{\gamma}$ has the property that for all $n$ such that $B_{n+1}$ through $B_{n+t}$ are constructed using the $\gamma$-Chacon-like operation and all $\frac{\tilde{h}_{n}}{M_{n}} -1 \leq \ell < \frac{\tilde{h}_{n+t}}{M_{n}}$,
\[
\mu(T^{\ell M_{n}}I_{n,0} \cap I_{n,0}) \leq  (1 - \kappa_{\gamma}) \mu(I_{n,0})
\]
for $\kappa_{\gamma}$ as in Proposition \ref{shifts}.
\end{lemma}
\begin{proof}
Since $M_{n}$ divides $|B_{n}1^{c_{n}}|$, treating $B_{n}1^{c_{n}}$ as $023\ldots P$ where $P = \frac{\tilde{h}_{n}}{M_{n}}$ and $1^{M_{n}}$ as $1$, the action of $T^{M_{n}}$ on $B_{n+1}$ through $B_{n+t}$ generates the same right-infinite word $W$ as the transformation in Proposition \ref{shifts}.

Suppose that $\mu(T^{\ell M_{n}}I_{n,0} \cap I_{n,0}) \geq (1 - \kappa_{\gamma} + \delta)\mu I_{n,0}$ for some $\delta > 0$.  Then the corresponding shift of $W$ has the property that $d(W_{[1,|W_{n}|]},(\sigma^{\ell}W)_{[1,|W_{n}|]}) \leq (\kappa_{\gamma} - \delta)|W_{n}|$ since it would force every block of $p$ levels to line up at least $(1 - \kappa_{\gamma} + \delta)$ of the time.
But Proposition \ref{shifts} states that $d(W_{[1,|W_{n}|]},(\sigma^{\ell}W)_{[1,|W_{n}|]}) \geq \kappa_{\gamma}|W_{n}|$ so this is a  contradiction.
\end{proof}

\begin{proposition}\label{pmix}
$T_{\gamma}$ is at least $\kappa_{\gamma}$-partially mixing where $\kappa_{\gamma} = \frac{2\gamma-1}{(2\gamma+2)(4\gamma+1)}$.
\end{proposition}
\begin{proof}
By construction, $r_{Q_{n}} = r_{Q_{n}+1} \to \infty$.  Also, since $|B_{t+1}| \geq 3|B_{t}|$ for all $t$,
\[
\frac{r_{Q_{n}}^{2}}{h_{Q_{n}}} = \frac{\tilde{h}_{S_{n-1}}^{2}}{h_{Q_{n}}} < \frac{1}{n} \frac{3^{Q_{n}}}{h_{Q_{n}}} \leq \frac{1}{n} \to 0
\]
and therefore the conditions in Proposition \ref{mix4} required on the $n$ for which the staircase operation is used are satisfied.

By Lemma \ref{Moverh}, $\frac{M_{n}}{h_{n-1}} \to 0$ along $n$ for which the $\gamma$-Chacon-like operation is used.  By construction, $c_{n}$ and $M_{n}$ satisfy the remaining conditions required by Proposition \ref{mix4} on the $n$ for which the Chacon-like operations are used.  The conditions required by Proposition \ref{pmix1} are also satisfied using the above and Lemma \ref{norigid}.

Let $D$ and $A$ be unions of levels in some fixed column $C_{m_{0}}$ and fix $\epsilon > 0$.  By Proposition \ref{mix4}, there exists $m_{1} > m_{0}$ such that for all $m \geq m_{1}$ for which $B_{m+1}$ is constructed using the staircase operation, $A$ and $D$ are both uniformly mixed within $\epsilon$ at stage $m$.  Therefore there is $n_{0}$ such that for all $n \geq n_{0}$, $A$ and $D$ are both uniformly mixed within $\epsilon$ at stages $Q_{n}$ and $Q_{n}+1$.

By Proposition \ref{pmix1}, since the $\gamma$-Chacon-like operation is used to construct $B_{S_{n+1}+1}$ through $B_{Q_{n+1}}$ and $S_{n+1} - 2 = Q_{n}$, then for all $(1 - \frac{1}{2}\epsilon)\tilde{h}_{Q_{n}} \leq t < (1 - \frac{1}{2}\epsilon)\tilde{h}_{Q_{n+1}}$,
\[
\mu(T^{t}D \cap A) \geq \kappa_{\gamma} \mu D \mu A - 5\epsilon
\]
As this holds for all $n \geq n_{0}$, then the above holds for all $t \geq \tilde{h}_{Q_{n_{0}}}$.  Since $\epsilon > 0$ was arbitrary, the claim is proved.
\end{proof}

\begin{proof}[Proof of Theorem \ref{examples}]
Let $\epsilon > 0$.  Choose $\gamma \geq 2$ such that $\frac{1}{\gamma+1} < \epsilon$.  By Proposition \ref{comp}, $T_{\gamma}$ satisfies $\liminf \nicefrac{p(q)}{q} \leq 2 + \frac{1}{\gamma+1} < 2 + \epsilon$.  By Proposition \ref{pmix}, $T_{\gamma}$ is at least $\kappa_{\gamma}$-partially mixing.
\end{proof}

\section{Lower bounds on word complexity of partially mixing rank-one subshifts}

\begin{definition}
A rank-one transformation is in \textbf{canonical form} if $B_{n}$ is a suffix of $B_{n+1}$ for all $n$.
\end{definition}

\begin{lemma}\label{canon}
Let $T$ be a rank-one transformation.  Then $T$ is measurably isomorphic to a rank-one transformation in canonical form  which generates the same language as $T$.
\end{lemma}
\begin{proof}
Since $B_{n+1}$ is a concatenation of $B_{n}$ and $1$, there exists a unique $c_{n}$ such that $B_{n+1}$ has $B_{n}1^{c_{n}}$ as a suffix.  Set $\tilde{B}_{n+1}$ to be the prefix of $B_{n+1}$ of length $|B_{n+1}| - c_{n}$.  Then $B_{n+1} = \tilde{B}_{n+1}1^{c_{n}}$.  Clearly $B_{n}$ and $\tilde{B}_{n}$ generate the same language.  The measure isomorphism is likewise immediate (and can be viewed as saying that we place zero spacers over the rightmost subcolumn but instead place that many extra above every subcolumn at the next stage).
\end{proof}

For the remainder of this section, all rank-one transformations are assumed to be in canonical form, without loss of generality.

\begin{definition}
A rank-one transformation has \textbf{bounded spacers} if there exists some $t \geq 0$ such that for every $n$ and $j$, $s_{n,j} \leq t$ and has \textbf{unbounded spacers} otherwise.
\end{definition}

\begin{lemma}\label{unbounded}
Let $T$ be a nontrivial rank-one subshift in canonical form that is lightly mixing.  Then $T$ has unbounded spacers.
\end{lemma}
\begin{proof}
Suppose that $T$ has bounded spacers.  Since $T$ is in canonical form, $T^{h_{n}}$ pushes every point in a level $I_{n,j}$ through exactly one spacer string (for points in the rightmost subcolumn, the spacers from a future stage).  Then $T^{h_{n}}I_{n,j} \subseteq \bigcup_{i=0}^{t} T^{-i}I_{n,j}$ for every height $h_{n}$.

Let $B$ be a union of levels in some fixed column such that every level in $B$ is at least $2t+2$ away from every other level in $B$.  Then $T^{h_{n}}B \cap T^{t+1}B = \varnothing$ since every point in $T^{h_{n}}B$ is within $t$ levels of a level in $B$ so none of those points can be exactly $t+1$ away as the levels in $B$ are all at least $2t+2$ apart.

This would contradict light mixing so $T$ must have unbounded spacers.  
\end{proof}

\begin{proof}[Proof of Theorem \ref{R1nolight}]
By Lemma \ref{canon}, we may assume $T$ is in canonical form.
By Lemma \ref{unbounded}, $T$ has unbounded spacers so $1^{\infty}$ is a left-infinite right-special word.  Therefore if for infinitely many $q$ we had $p(q+1) - p(q) = 1$ then the unique right-special word of length $q$ would be $1^{q}$.  The path from $1^{q}$ to itself starting with $1$ is labeled $1$, call the other path $u_{q}$.  Then $1 + |u_{q}| = p(q) + 1$ so $|u_{q}| = p(q)$.  Now $u_{q}$ must have $1^{q}$ as a suffix so, since the total measure of the spacers added is finite, it must be that $\nicefrac{q}{|u_{q}|} \to 0$.  But then $\nicefrac{p(q)}{q} \to \infty$.  So $p(q+1) - p(q) \geq 2$ for all sufficiently large $q$.  

Suppose that $\liminf \nicefrac{p(q)}{q} \leq 2$.  Fix $\epsilon > 0$.  By Theorem 2.2 \cite{boshernitzan}, there exist infinitely many $q$ such that $p(q+1) - p(q) \leq 2$ and $p(q) < (2 + \epsilon)q$.  For such $q$, by the above, $p(q+1) - p(q) = 2$ so there are two right-special words of length $q$.  Since $1^{q}$ must be one of them, call the other $w_{q}$.  Since $w_{q} \ne 1^{q}$, there is some $c_{q}, d_{q} \geq 1$ such that $w_{q} = 1^{d_{q}}z1^{c_{q}}$ for some word $z$ which starts and ends with $0$.

Since $1^{d_{q}}$ and $1^{c_{q}}$ correspond to spacers being added at some $n^{th}$ stage so $|w_{q}| \sim h_{n}$, and since $T$ is on a finite measure space, then $\nicefrac{c_{q}}{q} \to 0$ and $\nicefrac{d_{q}}{q} \to 0$.

Consider the Rauzy graph for such a $q$.  The path from $w_{q}$ to $1^{q}$ must have label $1^{q-d_{q}}$ because $w_{q}$ is followed by arbitrarily long strings of $1$s in the language.  The path from $1^{q}$ to $w_{q}$ must have label $z1^{c_{q}}$.  The path from $1^{q}$ to itself has label $1$.  Let $v_{q}$ be the label of the path from $w_{q}$ to itself.

Since the total length of the paths is $p(q) + 2$, then
\[
|v_{q}| + q - c_{q} + q - d_{q} + 1 = p(q) + 2 < (2+\epsilon)q + 2
\]
meaning that $|v_{q}| < \epsilon q + c_{q} + d_{q} + 1$.  

As this holds for all $\epsilon > 0$, then $\nicefrac{|v_{q}|}{q} \to 0$ over such $q$ (as $\nicefrac{c_{q}}{q} \to 0$ and $\nicefrac{d_{q}}{q} \to 0$).  By Lemma 5.6 \cite{creutzmixing}, $|v_{q}| < |w_{q}|$ implies $v_{q}$ is a root of $w_{q}$, meaning that $w_{q} = v_{q}^{\prime}v_{q}^{t_{q}}$ for some $t_{q} \geq 1$ and $v_{q}^{\prime}$ a suffix of $v_{q}$.  Since $\nicefrac{|v_{q}|}{q} = \nicefrac{|v_{q}|}{|w_{q}|} \to 0$, then $t_{q} \to \infty$.

Now $v_{q}$ is always followed by $v_{q}$ or by $1^{q - d_{q}}$ as those are the only paths in the Rauzy graph.  Note that if $\mu([1^{q-d_{q}}]) \geq \delta$ infinitely often then, as $T$ is rank-one, this would imply the measure of the spacers added at each level does not go to zero so the transformation would not be on a finite measure space.  So $v_{q}$ is followed by $1^{q-d_{q}}$ at most $\epsilon_{q}$ of the time for some $\epsilon_{q} \to 0$.

For a fixed $q_{0}$ and $q > q_{0}$, it then follows by the shift-invariance of the measure that for $x \in [v_{q_{0}}]$, we have $T^{|v_{q}|}x \in [v_{q_{0}}]$ for at least $(1 - \nicefrac{1}{t_{q}}- \epsilon_{q})\mu([v_{q_{0}}])$ since every occurrence of $v_{q_{0}}$ must be in some copy of $v_{q}$ and then $T^{|v_{q}|}$ will take that occurrence to the same occurrence in the following $v_{q}$ whenever $v_{q}$ is followed by $v_{q}$ (which must be at least $(1 - \nicefrac{1}{t_{q}}) - \epsilon_{q}$ of the time by shift-invariance and that $\mu([1^{q-d_{q}}]) \leq \epsilon_{q}$).  Therefore
\[
\mu(T^{|v_{q}|}[v_{q_{0}}] \cap [v_{q_{0}}]) \geq \left(1 - \frac{1}{t_{q}} - \epsilon_{q}\right)\mu([v_{q_{0}}])
\]
Since $\epsilon_{q} \to 0$ and $t_{q} \to \infty$, then
\[
\mu(T^{|v_{q}|}[v_{q_{0}}] \cap [v_{q_{0}}]) \to \mu([v_{q_{0}}])
\]
contradicting light mixing.
\end{proof}

\begin{remark}
In fact, the above proof shows that a rank-one transformation with $\liminf \nicefrac{p(q)}{q} = 2$ and unbounded spacers is rigid.  This follows since $\mu(\bigcup_{i=0}^{|v_{q_{0}}|} T^{i}[v_{q_{0}}]) \to 1$ since the spacer levels added have measure tending to zero.  However, we do not elaborate as we do not need this and there exist nonrigid rank-one transformations with $\limsup \frac{p(q)}{q} < 2$, e.g.~the examples in Section 4 \cite{creutzrankone} (which, of course, have bounded spacers)
\end{remark}

\begin{proof}[Proof of Theorem \ref{alphaoptimalR1}]
By Lemma \ref{canon}, we may assume $T$ is in canonical form.
We may and will assume $\liminf \frac{p(q)}{q} < 3 - 2\delta$ for some $\delta > 0$ since otherwise the claim is immediate.  Then, by Theorem 2.2 \cite{boshernitzan}, there exists infinitely many $q$ for which $p(q+1) - p(q) \leq 2$ and $p(q) \leq (3 - \delta)q$ and we may assume that the $\liminf$ is attained along this sequence of $q$'s.  Since $T$ must have unbounded spacers by Lemma \ref{unbounded}, $p(q+1) - p(q) \geq 2$ for all $q$.  So the $\liminf$ is attained along a sequence of $q$ such that $p(q+1) - p(q) = 2$ and $p(q) \leq (3 - \delta)q$.

For such $q$,
let $v_{q}$ be as in the proof of Theorem \ref{R1nolight}.  The paths in the Rauzy graph have total length
\[
|v_{q}| + 1 + q - c_{q} + q - d_{q} = p(q) + 2 \leq (3 - \delta)q + 2
\]
meaning that
\[
|v_{q}| \leq (1 - \delta)q + c_{q} + d_{q} + 1
\]
so again $v_{q}$ is a root of $w_{q}$ (as $\nicefrac{c_{q}}{q} \to 0$ and $\nicefrac{d_{q}}{q} \to 0$) so $w_{q} = v_{q}^{\prime}v_{q}^{t_{q}}$.

Let $C_{q} = \bigcup_{i=0}^{|v_{q}|-1} T^{i}[v_{q}]$ (note that we are not assuming the translates are disjoint).  Since $\mu([1^{q-d_{q}}]) \to 0$, it follows that $\mu(C_{q}) \to 1$.  As in the proof of Theorem \ref{R1nolight}, we have that for a fixed $q_{0}$ and any $q > q_{0}$ such that $p(q+1) - p(q) = 2$ and $p(q) \leq (3 - \delta)q$,
\[
\mu\left(T^{|v_{q}|}[v_{q_{0}}] \cap [v_{q_{0}}]\right) \geq \left(1 - \frac{1}{t_{q}} - \epsilon_{q}\right)\mu([v_{q_{0}}])
\]
and therefore
\[
\mu\left(T^{|v_{q}|}[v_{q_{0}}] \cap \left(C_{q_{0}} \setminus [v_{q_{0}}]\right)\right) < \left(\frac{1}{t_{q}} + \epsilon_{q}\right)\mu([v_{q_{0}}])
\]
Now $\mu(C_{q}) \to 1$ and $\mu([v_{q}]) \to 0$ (since $|v_{q}| \to \infty$ as otherwise the subshift is periodic), so $\mu(C_{q} \setminus [v_{q}]) \to 1$.

Fix $\epsilon > 0$ and choose $q_{0}$ in the set of possible $q$ such that $\mu(C_{q_{0}} \setminus [v_{q_{0}}]) > 1 - \epsilon$.  Then
\[
\mu\left(T^{|v_{q}|}[v_{q_{0}}] \cap \left(C_{q_{0}} \setminus [v_{q_{0}}]\right)\right) < \left(\frac{1}{t_{q}} + \epsilon_{q}\right)\mu([v_{q_{0}}])\frac{1}{1 - \epsilon}\mu\left(C_{q_{0}} \setminus [v_{q_{0}}]\right)
\] 
Since $T$ is $\alpha$-partially mixing,
\[
\mu\left(T^{|v_{q}|}[v_{q_{0}}] \cap \left(C_{q_{0}} \setminus [v_{q_{0}}]\right)\right)  \geq \alpha \mu([v_{q_{0}}])\mu\left(C_{q_{0}} \setminus [v_{q_{0}}]\right)
\]
meaning that, for $q$ large enough that $(1 - \epsilon)\alpha > \epsilon_{q}$,
\[
\alpha < \left(\frac{1}{t_{q}} + \epsilon_{q}\right)\frac{1}{1 - \epsilon}
\]
and therefore
\[
t_{q} < \frac{1}{(1 - \epsilon)\alpha - \epsilon_{q}}
\]
Since $p(q) \geq |v_{q}| + q - c_{q} + q - d_{q} + 1 - 2$,
\begin{align*}
\frac{p(q)}{q} &\geq \frac{|v_{q}|}{q} + 2 - \frac{d_{q} + c_{q} + 1}{q}
= \frac{|v_{q}|}{|w_{q}|} + 2 - \frac{d_{q} + c_{q} + 1}{q} \\
&> \frac{1}{t_{q}+1} + 2 - \frac{d_{q} + c_{q} + 1}{q}
> \frac{1}{\frac{1}{(1 - \epsilon)\alpha - \epsilon_{q}} + 1} + 2 - \frac{d_{q} + c_{q} + 1}{q}
\end{align*}
Taking $q \to \infty$ (so $\epsilon_{q} \to 0$ and $c_{q}/q \to 0$ and $d_{q} / q \to 0$) and $\epsilon \to 0$, then, along the sequence of $q$ such that $p(q+1) - p(q) = 2$,
\[
\liminf \frac{p(q)}{q} \geq 2 + \frac{\alpha}{\alpha + 1}
\]
so the claim follows as the $\liminf$ is attained along that sequence.
\end{proof}

\dbibliography{partialmixing}

\end{document}